\documentclass[a4paper,12pt]{amsart}
\usepackage{amssymb}
\author{Ivan V.~Arzhantsev}
\address{
Department of Higher Algebra\\
Faculty of Mechanics and Mathematics\\
Moscow State University\\
119992 Moscow, Russia}
\email{arjantse@mccme.ru}
\urladdr{http://mech.math.msu.su/department/algebra/staff/arzhan.htm}
\author{Dmitri A.~Timashev}
\address{
Department of Higher Algebra\\
Faculty of Mechanics and Mathematics\\
Moscow State University\\
119992 Moscow, Russia}
\email{timashev@mech.math.msu.su}
\urladdr{http://mech.math.msu.su/department/algebra/staff/timashev}
\thanks{Supported by RFBR grants 01--01--00756, 03--01--06252,
03--01--06253.}
\title[On canonical embeddings]{On the canonical embeddings of\\
certain homogeneous spaces}
\date{August 20, 2003}
\dedicatory{To A.~L.~Onishchik on his 70-th anniversary}

\newcommand{\kk}{\Bbbk}

\newcommand{\codim}{\mathop{\mathrm{codim}}\nolimits}
\newcommand{\trdeg}{\mathop{\mathrm{tr.deg}}}
\newcommand{\cchar}{\mathop{\mathrm{char}}}
\newcommand{\Aut}{\mathop{\mathrm{Aut}}}
\newcommand{\Spec}{\mathop{\mathrm{Spec}}}
\newcommand{\rk}{\mathop{\mathrm{rk}}}
\newcommand{\Ru}[1]{{#1}_{\text{\normalfont u}}}
\newcommand{\CE}{\mathop{\mathrm{CE}}}
\newcommand{\diag}{\mathop{\mathrm{diag}}}
\newcommand{\X}{\mathfrak{X}}
\newcommand{\dual}[1]{\check{#1}}
\newcommand{\ZZ}{\mathbb{Z}}
\newcommand{\QQ}{\mathbb{Q}}
\newcommand{\PP}{\mathbb{P}}
\newcommand{\n}{\mathfrak{n}}
\newcommand{\g}{\mathfrak{g}}
\newcommand{\h}{\mathfrak{h}}
\newcommand{\s}{\mathfrak{s}}
\newcommand{\z}{\mathfrak{z}}
\newcommand{\rr}{\mathfrak{r}}
\newcommand{\ssl}{\mathfrak{sl}}
\newcommand{\ggl}{\mathfrak{gl}}
\newcommand{\ad}{\mathop{\mathrm{ad}}}
\newcommand{\tr}{\mathop{\mathrm{tr}}}
\newcommand{\Hom}{\mathop{\mathrm{Hom}}}
\newcommand{\End}{\mathop{\mathrm{End}}}
\newcommand{\Mat}{\mathop{\mathrm{Mat}}}
\renewcommand{\Im}{\mathop{\mathrm{Im}}}
\newcommand{\Ker}{\mathop{\mathrm{Ker}}}
\newcommand{\conv}{\mathop{\mathrm{conv}}}
\newcommand{\1}{\boldsymbol{1}}
\newcommand{\md}{\mathop{\mathrm{mod}}\nolimits}
\newcommand{\res}[1]{\bar{#1}}
\newcommand{\resomega}{\res\omega}
\newcommand{\degs}{\dag}
\newcommand{\GL}{\mathrm{GL}}
\newcommand{\SL}{\mathrm{SL}}
\def\Sp{\mathrm{Sp}}
\newcommand{\Spin}{\mathrm{Spin}}
\newcommand{\Bb}{\mathbf{B}}
\newcommand{\Cc}{\mathbf{C}}
\newcommand{\Dd}{\mathbf{D}}
\newcommand{\Ee}{\mathbf{E}}
\newcommand{\Ff}{\mathbf{F}}
\newcommand{\Gg}{\mathbf{G}}

\newtheorem{theorem}{Theorem}
\newtheorem{corollary}{Corollary}
\newtheorem{proposition}{Proposition}
\newtheorem{lemma}{Lemma}
\newtheorem*{Claim}{Claim}
\theoremstyle{definition}
\newtheorem{definition}{Definition}
\newtheorem{example}{Example}
\newtheorem*{Example}{Example}
\newtheorem*{Examples}{Examples}
\theoremstyle{remark}
\newtheorem{remark}{Remark}

\begin{document}

\begin{abstract}
We study equivariant affine embeddings of homogeneous spaces and
their equivariant automorphisms. An example of a quasiaffine,
but not affine, homogeneous space with finitely many equivariant
automorphisms is presented. We prove the solvability of any
connected group of equivariant automorphisms for an affine
embedding with a fixed point and finitely many orbits. This is
applied to studying the orbital decomposition for algebraic
monoids and canonical embeddings of quasiaffine homogeneous
spaces, i.e., those affine embeddings associated with the
coordinate algebras of homogeneous spaces, provided the latter
algebras are finitely generated. We pay special attention to
the canonical embeddings of quotient spaces of reductive groups
modulo the unipotent radicals of parabolic subgroups. For these
varieties, we describe the orbital decomposition, compute the
modality of the group action, and find out which of them are
smooth.
We also describe minimal ambient modules for these canonical
embeddings provided that the acting group is simply connected.
\end{abstract}

\maketitle

\section{Introduction}

Let $G$ be a connected reductive algebraic group over an
algebraically closed field $\kk$ and $H$ be a closed subgroup of
$G$. It was proved by Y.~Matsushima~\cite{mat} and
A.~L.~Onishchik~\cite{on} that the homogeneous space $G/H$ is
affine if and only if $H$ is reductive. (For a simple proof,
see~\cite[\S2]{lu2}; a characteristic-free proof is given
in~\cite{rich}.) The subgroup $H$ is said to be \emph{observable}
in $G$ if the homogeneous space $G/H$ is a quasiaffine variety.
For a description of observable subgroups, see \cite{gr},
\cite{sukh}. In particular, any reductive subgroup is
observable.

Let us recall that $H$ is a \emph{Grosshans subgroup} in $G$ if
$H$ is observable and the algebra of regular functions
$\kk[G/H]$ is finitely generated.  This class of subgroups was
considered by F.~D.~Grosshans~\cite{gr1}, \cite{gr2}, \cite{gr}
in connection with the Hilbert 14-th problem. In particular, it
is proved in~\cite{gr2} that the unipotent radical $\Ru{P}$ of a
parabolic subgroup $P$ of $G$ is a Grosshans subgroup.

Let $H$ be a Grosshans subgroup in $G$.
The \emph{canonical embedding} of the homogeneous space $G/H$ is the
affine $G$-variety $\CE(G/H)=\Spec\kk[G/H]$ corresponding to the
affine algebra $\kk[G/H]$.  It is easy to see that $\CE(G/H)$ is
a normal affine variety with an open $G$-orbit isomorphic to
$G/H$, and the complement of the open orbit has codimension $\ge
2$. Moreover, these properties characterize $\CE(G/H)$ up to
$G$-equivariant isomorphism (for more details see~\cite{gr}).

The aim of this paper is to study the canonical embeddings of
the homogeneous spaces $G/\Ru{P}$. Such embeddings form a
remarkable class of affine quasi-homogeneous $G$-varieties. They
provide a geometric point of view at the properties of the
algebra $\kk[G]^{\Ru{P}}$.

We begin with the following general result on equivariant
automorphisms of an affine embedding $G/H\hookrightarrow X$: if
$X$ contains a $G$-fixed point and only finitely many
$G$-orbits, then the connected part of the group $\Aut_G(X)$ is
solvable (Theorem~\ref{teo1}). It is easy to prove that the
group of equivariant automorphisms of $\CE(G/H)$ is naturally
isomorphic to $N_G(H)/H$. We deduce that the number of
$G$-orbits in $\CE(G/\Ru{P})$ is infinite, except the trivial
cases (Proposition~\ref{pr2}).

A detailed description of $\CE(G/\Ru{P})$ is obtained in
Section~\ref{G/P_u} under the assumption $\cchar\kk=0$.

In fact, our approach works for any affine embedding
$G/\Ru{P}\hookrightarrow X$ with the maximal possible group of
$G$-equivariant automorphisms (equal to the Levi subgroup $L$
of~$P$). Such affine embeddings $X$ are classified by finitely
generated semigroups $S$ of $G$-dominant weights having the
property that all highest weights of tensor products of simple
$L$-modules with highest weights in $S$ belong to~$S$, too.
Furthermore, every choice of the generators
$\lambda_1,\dots,\lambda_m\in S$ gives rise to a natural
$G$-equivariant embedding $X\hookrightarrow\Hom(V^{\Ru{P}},V)$,
where $V$ is the sum of simple $G$-modules of highest weights
$\lambda_1,\dots,\lambda_m$, see Theorem~\ref{aff.emb}. The
convex cone $\Sigma^{+}$ spanned by $S$ is nothing else but the
dominant part of the cone $\Sigma$ spanned by the weight
polytope of $V^{\Ru{P}}$, see~\ref{affine.emb}.

We prove that the $(G\times L)$-orbits in $X$ are in bijection
with the faces of $\Sigma$ whose interiors contain dominant
weights, the orbit representatives being given by the projectors
onto the subspaces of $V^{\Ru{P}}$ spanned by eigenvectors of
eigenweights in a given face (Theorem~\ref{orb}). Also we
compute the stabilizers of these points in $G\times L$ and in
$G$, and the modality of the action $G:X$. Smooth embeddings are
classified by Theorem~\ref{smooth}.

These results are applied to canonical embeddings
$X=\CE(G/\Ru{P})$ as follows. The semigroup $S$ here consists of
all dominant weights, and $\Sigma$ is the span of the dominant
Weyl chamber by the Weyl group of~$L$. From Theorem~\ref{orb} we
deduce that $(G\times L)$-orbits in $X$ are in bijection with
the subdiagrams in the Dynkin diagram of $G$ such that no
connected component of such a subdiagram is contained in the
Dynkin diagram of $L$. In terms of these diagrams, we compute
the stabilizers and the modality of $G:X$, see
Corollary~\ref{orb.can}. From Theorem~\ref{smooth}, a
classification of smooth canonical embeddings stems
(Corollary~\ref{smooth(can)}).

The techniques used in the description of affine $(G\times
L)$-embeddings of $G/\Ru{P}$ are parallel to those
developed in \cite{grp.comp} for the study of equivariant
compactifications of reductive groups. This analogy becomes
more transparent in view of the bijection between these affine
embeddings $G/\Ru{P}\hookrightarrow X$ and algebraic monoids $M$
with the group of invertibles~$L$, given by $X=\Spec\kk[G\times^PM]$
(Proposition~\ref{emb<->mon}).

Finally, we describe the $G$-module structure on the tangent
space of $\CE(G/\Ru{P})$ at the $G$-fixed point, assuming that
$G$ is simply connected simple. This space is obtained from
$\bigoplus_i\Hom(V_i^{\Ru{P}},V_i)$, where $V_i$ are the
fundamental simple $G$-modules, by removing certain summands
according to an explicit algorithm, see Theorem~\ref{tangent}.
The tangent space at the fixed point is at the same time the
minimal ambient $G$-module for $\CE(G/\Ru{P})$.
Its dual space is canonically isomorphic to the linear span of a
minimal homogeneous generating set for the algebra $\kk[G]^{\Ru{P}}$,
which is positively graded unless $P=G$.

Aside from the main subject of the paper, in
subsection~\ref{Losev} we provide an example of an observable
non-reductive subgroup $H\subset G$ such that the group
$N_G(H)/H$, and hence $\Aut_G(X)$ for any embedding
$G/H\hookrightarrow X$, is finite. This example answers the
conjecture in \cite{at} in the negative. We are grateful to
I.~V.~Losev for this example.

\subsection*{Acknowledgement}

This work was started during the stay of both authors at
Institut Fourier in spring~2003. We would like to thank this
institution for hospitality and Prof.~M.~Brion for invitation.


\section{Equivariant automorphisms}

Let $G/H$ be a homogeneous space. By $N_G(H)$ denote the
normalizer of $H$ in $G$. The group $\Aut_G(G/H)$ of
$G$-equivariant automorphisms of $G/H$ is isomorphic to
$N_G(H)/H$, where $nH$ acts on $G/H$ by $nH*gH=gn^{-1}H$,
$\forall n\in N,\ g\in G$.

Recall that an \emph{affine embedding} of a homogeneous space $G/H$ is an
affine $G$-variety $X$ containing a point $x\in X$ such that the orbit $Gx$ is
dense in $X$ and the orbit morphism $G\to Gx,\ g\mapsto gx$
induces an isomorphism between $G/H$ and $Gx$. In this situation
we use the notation $G/H\hookrightarrow X$. The embedding is
said to be \emph{trivial} if $Gx=X$.

\sloppy

\subsection{Automorphisms}\label{Aut_G}

For an embedding $G/H\hookrightarrow X$,
the group $\Aut_G(X)$ preserves the open orbit and may be considered
as a (closed) subgroup of $N_G(H)/H$.

\fussy

It is natural to ask which subgroups of $N_G(H)/H$ can be realized as
$\Aut_G(X)$, where $X$ is as above. Let us list some results in
this direction, assuming $\cchar\kk=0$:

\begin{enumerate}
\renewcommand{\theenumi}{\textup{(\arabic{enumi})}}
\renewcommand{\labelenumi}{\theenumi}
\sloppy

\item if $G/H$ is a spherical homogeneous space, then $\Aut_G(X)=N_G(H)/H$
for any affine embedding $G/H\hookrightarrow X$, see
e.g.~\cite[\S5]{at};

\item\label{Aut(SL(2)-emb)} if $G=\SL(2)$, $H=\{e\}$, then for
any non-trivial normal affine embedding the group $\Aut_G(X)$ is
a Borel subgroup in $N_G(H)/H\cong\SL(2)$~\cite{pop};

\item if $H$ is reductive, then the following conditions are
equivalent (cf.~\cite[Prop.\,2]{at}):

\begin{itemize}

\item for any non-trivial affine
embedding $G/H\hookrightarrow X$ one has $\dim\Aut_G(X)<\dim N_G(H)/H$;

\item $N_G(H)/H$ is a semisimple group.

\end{itemize}

\fussy

Indeed, let $L^0$ denote the connected component of unit in an
algebraic group~$L$. An affine embedding $G/H\hookrightarrow X$
such that $\dim\Aut_G(X)=\dim N_G(H)/H$ may be regarded as a
$\widehat{G}$-equi\-vari\-ant embedding of
$\widehat{G}/\widehat{H}$, where
$\widehat{G}=G\times(N_G(H)/H)^0$ and $\widehat{H}=\{\,(n,nH)\mid
n\in N_G(H),\ nH\in(N_G(H)/H)^0\,\}$. If $N_G(H)/H$ is not
semisimple, then $(N_G(H)/H)^0$ contains a central
one-dimensional torus~$S$, whence
$N_{\widehat{G}}(\widehat{H})/\widehat{H}\supseteq\widehat{S}=
\{e\}\times S$. Let $\widehat{N}\subseteq
N_{\widehat{G}}(\widehat{H})$ be the extension of $\widehat{S}$
by~$\widehat{H}$. Then there exists a non-trivial embedding
$\widehat{G}/\widehat{H}\hookrightarrow
X=\widehat{G}\times^{\widehat{N}}\mathbb{A}^1$, where the
quotient torus $\widehat{S}=\widehat{N}/\widehat{H}$ acts on
$\mathbb{A}^1$ by homotheties. This proves the direct
implication. The converse implication stems from Luna's
theorem~\cite{lu1}, since
$N_{\widehat{G}}(\widehat{H})/\widehat{H}$ is finite if
$N_{G}(H)/H$ is semisimple.

\end{enumerate}

The main result of this section may be considered as a partial generalization
of item~\ref{Aut(SL(2)-emb)}.

\begin{theorem}\label{teo1}
Let $G/H\hookrightarrow X$ be an affine embedding with a finite number of
$G$-orbits and with a $G$-fixed point. Then the group $\Aut_G(X)^0$
is solvable.
\end{theorem}

We begin the proof with the following

\begin{lemma}\label{lem1}
Let $X$ be an affine variety with an action of a connected
semisimple group $S$. Suppose that there is a point $x\in X$
and a one-parameter subgroup $\gamma:\kk^{\times}\to S$ such that
$\lim_{t\to 0}\delta(t)x$ exists in $X$ for any subgroup $\delta$
conjugate to~$\gamma$. Then $x$ is a $\gamma(\kk^{\times})$-fixed
point.
\end{lemma}

\begin{proof}
Let $T$ be a maximal torus in $S$ containing
$\gamma(\kk^{\times})$. It is known (for example,
see~\cite{pv2}) that $X$ can be realized as a closed $S$-stable
subvariety in $V$ for a suitable $S$-module $V$.  Let
$x=x_{\lambda_1}+\dots+x_{\lambda_n}$ be the weight
decomposition (with respect to~$T$) of $x$ with weights
$\lambda_1,\dots,\lambda_n$.  One-parameter subgroups of $T$
form the lattice $\X_{*}(T)$ dual to the character lattice
$\X(T)$. The existence of $\lim_{t\to 0}\gamma(t)x$ in $X$ means
that all pairings $\langle\gamma,\lambda_i\rangle$ are
non-negative.  Let $\gamma_1,\dots,\gamma_m$ be all the
translates of $\gamma$ under the action of the Weyl group
$W=N_S(T)/T$. By assumption,
$\langle\gamma_j,\lambda_i\rangle\ge 0$ for any $i=1,\dots,n$,
$j=1,\dots,m$, hence
$\langle\gamma_1+\dots+\gamma_m,\lambda_i\rangle\ge 0$.  Since
$\gamma_1+\dots+\gamma_m=0$, one has
$\langle\gamma_j,\lambda_i\rangle=0$ for any $i,\ j$. This shows
that the points~$x_{\lambda_i}$ (and~$x$) are
$\gamma(\kk^{\times})$-fixed.
\end{proof}

The following
proposition is a generalization of~\cite[Thm.\,4.3]{gr2}.

\begin{proposition}\label{prop1}
Suppose that $G/H\hookrightarrow X$ is an affine embedding with
a non-trivial $G$-equivariant action of a connected semisimple
group~$S$. Then the orbit $S*x$ is closed in $X$, $\forall x\in
G/H$.
\end{proposition}

\begin{proof}
We may assume $x=eH$. If $S*x$ is not closed, then,
by~\cite[Thm.\,1.4]{ke}, there is a one-parameter subgroup
$\gamma:\kk^{\times}\to S$ such that the limit
$$
  \lim_{t\to 0}\gamma(t)*x
$$
exists in $X$ and does not belong to $S*x$. Replacing $S$ by a
finite cover, we may assume that $S$ embeds in $N_G(H)$ (and
thus in~$G$) with a finite intersection with~$H$. By the
definition of $*$-action, one has
$\gamma(t)*x=\gamma(t^{-1})x$.
For any $s\in S$ the limit
$$
\lim_{t\to 0} (s\gamma(t))*x=\lim_{t\to 0}\gamma(t^{-1})s^{-1}x
$$
exists. Hence
$\lim_{t\to 0}s\gamma(t^{-1})s^{-1}x$ exists, too.
This shows that
for any one-parameter subgroup $\delta$ of $S$,
conjugate to $-\gamma$,
$\lim_{t\to 0}\delta(t)x$ exists in~$X$.
Lemma~\ref{lem1} implies that $x=\lim_{t\to 0}\gamma(t)*x$,
and this contradiction proves
Proposition~\ref{prop1}.
\end{proof}

\begin{proof}[Proof of the theorem]
Suppose that $\Aut_G(X)^0$ is not solvable. Then there is
a connected semisimple group $S$ acting on $X$ $G$-equivariantly.
By Proposition~\ref{prop1}, any $(S,*)$-orbit in the open
$G$-orbit of $X$ is closed in~$X$. In particular, the
$(S,*)$-action on $X$ is stable.

Let $X_1$ be the closure of a $G$-orbit in~$X$. Since
$G$ has a finite number of orbits in $X$, the variety $X_1$ is
$(S,*)$-stable. Applying the above arguments to~$X_1$, we
show that any $(S,*)$-orbit in $X$ is closed. But in this case
all $(S,*)$-orbits have the same dimension $\dim S$. On the
other hand, a $G$-fixed point is an $(S,*)$-orbit, a
contradiction.
\end{proof}

\begin{corollary}[of the proof]\label{c1}
Let $X$ be an affine $G$-variety with an open $G$-orbit. Suppose that
\begin{enumerate}
\renewcommand{\theenumi}{\textup{(\arabic{enumi})}}
\renewcommand{\labelenumi}{\theenumi}

\item a semisimple group $S$ acts on $X$ effectively and $G$-equivariantly;

\item\label{dim(closed)} the dimension of a closed $G$-orbit in
$X$ is less than $\dim S$.

\end{enumerate}

Then the number of $G$-orbits in $X$ is infinite.
\end{corollary}

\begin{remark}
Condition~\ref{dim(closed)} is essential. Indeed, let $H$ be
a one-di\-men\-sional unipotent root subgroup of $G=\SL(n)$. Then
$X=\CE(G/H)\cong\SL(n)\times^{\SL(2)}\mathbb{A}^2$, where $\SL(2)$
embeds in $\SL(n)$ as the standard 3-di\-men\-sional simple
subgroup containing~$H$, has two orbits, and $S=\SL(n-2)\subset
N_G(H)/H$.
\end{remark}

\subsection{Example of Losev}\label{Losev}

In many cases, Theorem~\ref{teo1} may be used to show that the group
$\Aut_G(X)$ cannot be very big.
%
On the other hand, there exist quasiaffine homogeneous spaces
$G/H$ such that $N_G(H)/H$ is finite and therefore $\Aut_G(X)$
is finite for every embedding $G/H\hookrightarrow X$. Such
examples with affine $G/H$ are well known: for instance, this is
the case if $H$ is a reductive subgroup containing a maximal
torus of~$G$. In fact, if $H$ is reductive and $N_G(H)/H$
finite, then there exists only a trivial affine embedding
$X=G/H$ \cite{lu1}.

It was conjectured in \cite{at} that $N_G(H)/H$ is infinite
whenever $H\subseteq G$ is observable, but not reductive.
However in 2003 I.~V.~Losev found a counterexample, which we
reproduce here with his kind permission.

Assume $\cchar\kk=0$. A desired subgroup $H\subseteq G$ is
sought in the form $H=R\leftthreetimes S$, where $R$ is
non-trivial unipotent and $S$ connected semisimple. Such a
subgroup $H$ has no non-trivial characters, whence $H$ is
observable \cite[1.5]{gr}.

We denote the Lie algebras of algebraic groups by the respective
lowercase Gothic letters. In order to obtain that $N_G(H)/H$ is
finite, it suffices to construct a unipotent subalgebra
$\rr\subset\g$ and a semisimple subalgebra
$\s\subset\n_{\g}(\rr)$ such that $\n_{\g}(\h)=\h$, where
$\h=\rr+\s$.

Note that if $\n_{\g}(\h)\ne\h$, then there is an element
$x\in\n_{\g}(\h)\setminus\h$ such that $[x,\s]\subseteq\s$. Any
derivation of $\s$ is inner, hence we may suppose that
$x\in\z_{\g}(\s)$. Thus for our purposes it suffices to
construct $\rr,\s$ such that $[x,\rr]\nsubseteq\rr$, $\forall
x\in\z_{\g}(\s)\setminus\{0\}$.

Take a simple Lie algebra~$\s$. We fix Cartan and Borel
subalgebras in $\s$ and denote by $V(\lambda)$ a simple
$\s$-module of highest weight~$\lambda$. The highest weight of
$V(\lambda)^{*}$ is denoted by~$\lambda^{*}$. Take three
distinct dominant weights $\lambda,\mu,\nu$ and put
$V=V(\lambda)\oplus V(\mu)\oplus V(\nu)$. We have an embedding
$\s\hookrightarrow\g=\ssl(V)$. The $\ad(\s)$-module structure of
$\ggl(V)$ is represented at the following picture:
\begin{equation*}
\ggl(V)=\left(
\begin{array}{c|c|c}
\ggl(V(\lambda)) & V(\lambda)\otimes V(\mu^*) & V(\lambda)\otimes V(\nu^*) \\
\hline
V(\lambda^*)\otimes V(\mu) & \ggl(V(\mu)) & V(\mu)\otimes V(\nu^*) \\
\hline
V(\lambda^*)\otimes V(\nu) & V(\mu^*)\otimes V(\nu) & \ggl(V(\nu)) \\
\end{array}\right)
\end{equation*}

Suppose that there exists a dominant weight $\rho\ne0$
satisfying the following conditions:
\begin{enumerate}
\item[($\rho1$)] $V(\rho)$ embeds into $V(\lambda)\otimes
V(\mu^*)$, $V(\lambda)\otimes V(\nu^*)$, $V(\mu)\otimes
V(\nu^*)$ as an $\s$-submodule;
\item[($\rho2$)] $V(\rho)$ does not embed into $V(\rho)\otimes
V(\rho)$.
\end{enumerate}
Consider the diagonal embedding of $V(\rho)$ in the direct sum
of $V(\lambda)\otimes V(\mu^*)$, $V(\lambda)\otimes V(\nu^*)$,
$V(\mu)\otimes V(\nu^*)$. Its image $\rr_1$ may be naturally
identified with a subspace in $\ssl(V)$ such that
$\rr_1\cap[\rr_1,\rr_1]=0$. Now
$\rr=\rr_1+[\rr_1,\rr_1]\subset\ssl(V)$ is a unipotent
subalgebra and $[\s,\rr]\subseteq\rr$, $[\s,\rr_1]=\rr_1$.

\begin{Claim}
We have $\z_{\g}(\s)\cap\n_{\g}(\rr)=0$. Thus
$\h=\rr+\s\subset\g$ is the desired subalgebra.
\end{Claim}
\begin{proof}
Clearly, $\z_{\g}(\s)$~is the two-dimensional diagonal toric
subalgebra of traceless block-scalar matrices. Any element
$x\in\z_{\g}(\s)$ is a diagonal matrix with
$x|_{V(\lambda)}=x_1\cdot\1_{V(\lambda)}$,
$x|_{V(\mu)}=x_2\cdot\1_{V(\mu)}$,
$x|_{V(\nu)}=x_3\cdot\1_{V(\nu)}$. Then $\ad x$ acts on
$V(\lambda)\otimes V(\mu^*)$, $V(\lambda)\otimes V(\nu^*)$,
$V(\mu)\otimes V(\nu^*)$ by constants $x_1-x_2$, $x_1-x_3$,
$x_2-x_3$, respectively.

By the condition ($\rho2$), if $[x,\rr]\subseteq\rr$, then
$[x,\rr_1]\subseteq\rr_1$, hence $x_1-x_2=x_1-x_3=x_2-x_3$,
i.e., $x_1=x_2=x_3$. The condition $\tr x=0$ implies $x=0$.
\end{proof}

So it remains to find dominant weights $\lambda,\mu,\nu,\rho$
satisfying the conditions ($\rho1$) and~($\rho2$).

It is known \cite{PRV} that the multiplicity of a simple
$\s$-submodule $V(\rho)$ in $V(\lambda)\otimes V(\mu)$ is equal
to
\begin{equation*}
\dim\{v\in V(\rho)_{\lambda-\mu^*} \mid e_i^{\mu^*_i+1}v=0,\
\forall i\}
\end{equation*}
Here $V(\rho)_{\lambda-\mu^*}$ is the weight subspace in
$V(\rho)$ of eigenweight $\lambda-\mu^*$, $\mu^*_i$~is the
numerical label of $\mu^*$ at the simple root~$\alpha_i$, and
$e_i$ is a non-zero element in~$\s_{\alpha_i}$.

Note that if $\mu^*_i\ge\rho'_i$, where $\rho'$ runs over the
orbit of $\rho$ under the Weyl group, then the multiplicity of
$V(\rho)$ in $V(\lambda)\otimes V(\mu)$ equals $\dim
V(\rho)_{\lambda-\mu^*}$. Thus if we add a weight with
sufficiently big numerical labels to $\lambda$, $\mu$, and
$\nu$, then the condition ($\rho1$) can be reformulated as
follows: $\lambda-\mu$, $\lambda-\nu$, $\mu-\nu$ occur among the
weights of~$V(\rho)$. In order to verify~($\rho2$), it suffices
to show that $\rho-\rho^*$ does not occur among the weights
of~$V(\rho)$.

\begin{Example}
Take $\s=\ssl(n)$, $n>2$. Let $\omega_1,\dots,\omega_{n-1}$
denote the fundamental weights and
$\varepsilon_1,\dots,\varepsilon_n$ the weights of the
tautological representation~$V(\omega_1)$.  Take
$\rho=2n\omega_1$ and $\lambda,\mu,\nu$ such that
$\lambda-\mu=\mu-\nu=n\omega_1$, $\lambda-\nu=2n\omega_1$. In
order to have numerical labels big enough, it suffices to take
$\nu=2n(\omega_1+\dots+\omega_{n-1})$.

As $\lambda-\mu$, $\lambda-\nu$, $\mu-\nu$ are vectors in the
weight polytope of $V(\rho)$ congruent to $\rho$ modulo the root
lattice, they occur among the weights of~$V(\rho)$. Finally,
$\rho-\rho^*=2n(\varepsilon_1+\varepsilon_n)$ is not a weight
of~$V(\rho)$, because the dominant weight $2n\omega_2$ in its
orbit under the Weyl group is not a weight of~$V(2n\omega_1)$.
\end{Example}

\subsection{Canonical embeddings}\label{can.emb}

Now we are going to apply the obtained results to the study of
the canonical embedding $\CE(G/\Ru{P})$. Fix a pair $T\subset B$, where $B$
is a Borel subgroup in $G$ and $T$ is a maximal torus. We shall consider
a parabolic subgroup $P\supseteq B$.

\begin{remark}\label{G->simple}
The commutator subgroup $G'\subseteq G$ is the maximal
semisimple subgroup in~$G$. It is easy to see that
$\CE(G/\Ru{P})=G\times^{G'}\CE(G'/\Ru{P})$ is the homogeneous
fibration over $G/G'$ with fiber $\CE(G'/\Ru{P})$. Thus without
loss of generality we may assume $G$ to be semisimple.

Furthermore,
$\CE(G/\Ru{P})=\CE(\widetilde{G}/\Ru{\widetilde{P}})/\widetilde{Z}$,
where $\widetilde{G}\to G=\widetilde{G}/\widetilde{Z}$ is the
simply connected covering, and
$\widetilde{P}\subseteq\widetilde{G}$ is the preimage of~$P$.
Passing to the quotient modulo a finite central subgroup
preserves many features of $\CE(G/\Ru{P})$ (for instance, the
orbit dimensions, the modality of the $G$-action, normality,
etc). Therefore we may assume in many questions that $G$ is
simply connected.

Then $G=\prod_iG_i$, $P=\prod_iP_i$, where $G_i$ are simple
factors. It follows that
$\CE(G/\Ru{P})=\prod_i\CE\bigl(G_i/\Ru{(P_i)}\bigr)$, and we may assume
$G$ to be a simply connected simple algebraic group.
\end{remark}

The following proposition gives a partial answer to a question
posed in~\cite{ar2}.

\begin{proposition}\label{pr2}
The number of $G$-orbits in $\CE(G/\Ru{P})$ is finite if and only if either
$P\cap G_i=G_i$ or $P\cap G_i=B\cap G_i$ for each simple factor
$G_i\subseteq G$.
\end{proposition}

\begin{proof}
We may assume by Remark~\ref{G->simple} that $G$ is simple.
If $P=G$, then $\Ru{P}=\{e\}$ and $\CE(G/\Ru{P})=G$. If
$P=B$, then $\Ru{P}=U$ is a maximal unipotent subgroup of $G$, the variety
$\CE(G/\Ru{P})$ is spherical \cite[2.1]{var.sph}, and any
spherical variety contains a finite number of $G$-orbits
[loc.~cit.].

To prove the converse implication, let us fix some
notation. Let $L$ be the Levi subgroup of $P$ containing~$T$,
$Z$~the center of~$L$, and $S$ the maximal semisimple subgroup
of~$L$.

\begin{lemma}\label{fixed pt}
If $P\ne G$, then $\CE(G/\Ru{P})$ contains a $G$-fixed point.
\end{lemma}

\begin{proof}[Proof of the lemma]
It is possible to find a one-parameter subgroup
$\gamma:\kk^{\times}\to Z$, such that the pairing of $\gamma$
with any non-zero dominant weight is positive (see e.g.\
Remark~\ref{grading} below). The 1-torus $\gamma(\kk^{\times})$,
considered as a subgroup of $\Aut_G(G/\Ru{P})$, defines a
$G$-invariant grading on $\kk[G/\Ru{P}]$.

The homogeneous
subalgebra $\kk[G/U]\subseteq\kk[G/\Ru{P}]$ has a multi-grading
$\kk[G/U]=\bigoplus_{\lambda}E(\lambda)$ by eigenspaces of the
$T$-action from the right, which are called \emph{dual Weyl
modules}, so that $\deg E(\lambda)=
\langle\gamma,\lambda\rangle$. It is known that $E(\lambda)\ne0$
iff $\lambda$ is dominant, and $\dim E(\lambda)=1\iff\lambda=0$
(see e.g.~\cite[\S12]{gr}). Hence the grading is non-negative on
$\kk[G/U]$ and only constant functions have degree~$0$.

The same
is true for the algebra $A=\kk[S*M]$, where $M$ is any
homogeneous generating set of~$\kk[G/U]$, since the
${*}$-actions of $S$ and $\gamma$ commute. However
$\kk[G/\Ru{P}]$ is the integral closure of~$A$ \cite[5.4]{gr2}
($A=\kk[G/\Ru{P}]$ if $\cchar\kk=0$, cf.~\ref{k[G/P_u]}). It
easily follows that the grading is non-negative on $\kk[G/\Ru{P}]$,
and the positive part of this grading is a $G$-stable maximal
ideal in $\kk[G/\Ru{P}]$.
\end{proof}

Note that the group $N_G(\Ru{P})/\Ru{P}$ is isomorphic to $L$. Hence if
$P\ne B$, then $S$ acts on $\CE(G/\Ru{P})$
effectively and $G$-equivariantly. By Theorem~\ref{teo1}, the
number of $G$-orbits in $\CE(G/\Ru{P})$ is infinite.
\end{proof}

\subsection{Reductive monoids}\label{red.mon}

Finally, consider an application of Theorem~\ref{teo1} to another remarkable
class of affine embeddings.

\begin{proposition}\label{pr3}
Let $M$ be an algebraic monoid with zero such that its group of
invertible elements $G(M)$ is reductive. Then
the number of left $G(M)$-cosets in $M$ is finite if and
only if $M$ is commutative.
\end{proposition}

\begin{proof}
It is known that $M$ is an affine variety \cite{mon} and $G(M)$
is open in $M$~\cite{vi}. If $M$ is commutative, then $G(M)$ is
commutative, hence $M$ is a toric variety and the number of
$G(M)$-orbits in $M$ is finite.

Otherwise, $G(M)$ contains a semisimple subgroup $S$, and the action of $S$ on
$M$ by right multiplication is $G(M)$-equivariant. The zero element is a
$G(M)$-fixed point, and we conclude by Theorem~\ref{teo1}.
\end{proof}


\section{The canonical embedding of $G/\Ru{P}$}
\label{G/P_u}

In this section we obtain a detailed description of
$\CE(G/\Ru{P})$, assuming $\cchar\kk=0$. Our basic idea is to
consider $G/\Ru{P}$ as a homogeneous space under $G\times L$, in
the notation of the previous section. The action is defined by
$(g,h)x\Ru{P}=gxh^{-1}\Ru{P}$, $\forall g,x\in G$, $h\in L$.

It is clear that $X=\CE(G/\Ru{P})$ is a $(G\times
L)$-equivariant affine embedding of $G/\Ru{P}$. More generally,
we shall describe the structure of an arbitrary affine $(G\times
L)$-embedding of $G/\Ru{P}$ and deduce results concerning the
canonical embedding as a particular case.

\subsection{The coordinate algebra}\label{k[G/P_u]}

One easily sees from the Bruhat decomposition that a Borel
subgroup of $G\times L$ has an open orbit in $G/\Ru{P}$, i.e.,
$G/\Ru{P}$ is a spherical homogeneous $(G\times L)$-space.
Alternatively, one can deduce that $G/\Ru{P}$ is spherical from
the multiplicity-free property for the isotypic decomposition of
$\kk[G/\Ru{P}]$, which we are going to describe.

Let $V(\lambda)=V_G(\lambda)$ denote the simple $G$-module of
highest weight $\lambda$ w.r.t.~$B$, where $\lambda$ is any
vector from the semigroup $\X^{+}=\X^{+}_G$ of $B$-dominant
weights. The set of positive/negative roots w.r.t.~$B$ is
denoted by $\Delta^{\pm}=\Delta_G^{\pm}$, and
$\Pi=\Pi_G\subseteq\Delta_G^{+}$ is the set of simple roots. The
respective sets of coroots are denoted by
$\dual\Delta^{\pm}$, $\dual\Pi$. Let $C=C_G$ and $\dual{C}$ be
the dominant Weyl chambers in $\X(T)\otimes\QQ$ and
$\X_{*}(T)\otimes\QQ$, respectively. Recall that $C,\dual{C}$
are fundamental chambers of the Weyl group $W=W_G=N_G(T)/T$.

The group $G$ itself can be considered as a homogeneous space
$(G\times G)/\diag{G}$, where the left and right copies of $G$
act by left and right translations, respectively. The $(G\times
G)$-isotypic decomposition of $\kk[G]$ is well known:

\begin{proposition}[{\cite[II.3.1, Satz~3]{inv}}]\label{k[G]}
$\kk[G]=\bigoplus_{\lambda\in\X^{+}}\kk[G]_{(\lambda)}$, where
$\kk[G]_{(\lambda)}\cong V(\lambda)^{*}\otimes V(\lambda)$ is the
linear span of the matrix elements of the representation
$G:V(\lambda)$.
\end{proposition}

Now the isotypic decomposition of $\kk[G/\Ru{P}]$ is provided by
passing to $\Ru{P}$-invariants from the right in the r.h.s.\ of
the above decomposition. Note that $V(\lambda)^{\Ru{P}}\cong
V_L(\lambda)$ is isomorphic to the simple $L$-module of highest
weight~$\lambda$.

\begin{proposition}
There is a $(G\times L)$-module decomposition
\begin{equation*}
\kk[G/\Ru{P}]=\bigoplus_{\lambda\in\X^{+}}\kk[G/\Ru{P}]_{(\lambda)}
\end{equation*}
where $\kk[G/\Ru{P}]_{(\lambda)}\cong V(\lambda)^{*}\otimes
V_L(\lambda)$ is the linear span of the matrix elements of the
linear maps $V(\lambda)^{\Ru{P}}\to V(\lambda)$ induced by $g\in
G$, considered as regular functions on~$G/\Ru{P}$.
\end{proposition}

\begin{proof}
In view of Proposition~\ref{k[G]} and the above remark, it
suffices to note that the space of matrix elements of linear
maps $V(\lambda)^{\Ru{P}}\to V(\lambda)$ equals
$\left(\kk[G]_{(\lambda)}\right)^{\Ru{P}}$ (invariants under
right translations).
\end{proof}

Next we describe the multiplicative structure of $\kk[G/\Ru{P}]$.

\begin{proposition}\label{tails}
There is a decomposition
\begin{equation*}
\kk[G/\Ru{P}]_{(\lambda)}\cdot\kk[G/\Ru{P}]_{(\mu)}=
\kk[G/\Ru{P}]_{(\lambda+\mu)}
\oplus\bigoplus_i\kk[G/\Ru{P}]_{(\lambda+\mu-\beta_i)}
\end{equation*}
$\forall\lambda,\mu\in\X^{+}$, where $\lambda+\mu-\beta_i$ runs
over the highest weights of all ``lower'' irreducible components
in the $L$-module decomposition $V_L(\lambda)\otimes V_L(\mu)=
V_L(\lambda+\mu)\oplus\cdots$, so that $\beta_i\in\ZZ_{+}\Pi_L$.
\end{proposition}

\begin{proof}
$\kk[G/\Ru{P}]_{(\lambda)}\cdot\kk[G/\Ru{P}]_{(\mu)}$ is spanned
by the products of matrix elements of linear maps
$V(\lambda)^{\Ru{P}}\to V(\lambda)$ and $V(\mu)^{\Ru{P}}\to
V(\mu)$ induced by $g\in G$, i.e., by matrix elements of
$V(\lambda)^{\Ru{P}}\otimes V(\mu)^{\Ru{P}}\to V(\lambda)\otimes
V(\mu)$. But $V(\lambda)^{\Ru{P}}\otimes V(\mu)^{\Ru{P}}\cong
V_L(\lambda)\otimes V_L(\mu)$, and each $L$-highest weight
vector occurring in the l.h.s.\ is a $G$-highest weight vector
at the same time, because it is fixed by~$\Ru{P}$. It generates
a simple $G$-submodule $V(\lambda+\mu-\beta)\subseteq
V(\lambda)\otimes V(\mu)$ and a simple $L$-submodule
$V_L(\lambda+\mu-\beta)=V(\lambda+\mu-\beta)^{\Ru{P}}$, where
$\beta\in\ZZ_{+}\Pi_L$. The latter $L$-submodule is mapped to
$V(\lambda+\mu-\beta)$ by $g\in G$. Therefore the above space of
matrix elements for tensor products is spanned by all the
$\kk[G/\Ru{P}]_{(\lambda+\mu-\beta)}$.
\end{proof}

\begin{remark}\label{grading}
Let $\gamma$ be any vector in the interior of the cone
$\dual{C}\cap\Pi_L^{\perp}\cap\langle\dual\Pi\rangle$. For
instance, one may take $\gamma=\dual\rho_L$, the sum of the
fundamental coweights corresponding to simple roots from
$\Pi\setminus\Pi_L$, or $\gamma=\dual\rho_G-\dual\rho_L=
\frac12\sum_{\alpha\in\Delta_G^{+}\setminus\Delta_L^{+}}\dual\alpha$.
(Here
$\dual\rho_G=\frac12\sum_{\alpha\in\Delta_G^{+}}\dual\alpha$ is
half the sum of positive coroots, or equivalently, the sum of
fundamental coweights of a reductive group~$G$.) Then
$\langle\gamma,C\rangle\ge0$, and the inequality is strict on
$C\setminus\{0\}$ if $G$ is simple. (This is because
$C\setminus\{0\}$ is contained in the interior of the cone
$\QQ_{+}\Pi$ dual to $\dual{C}$, for indecomposable root
systems.)

Replacing $\gamma$ by a multiple, we may assume that
$\gamma\in\X_{*}(T)$ defines a one-parameter subgroup
$\kk^{\times}\to Z$. This 1-torus defines an invariant
non-negative algebra grading of $\kk[G/\Ru{P}]$ via the action
by right translations of an argument, so that
$\deg\kk[G/\Ru{P}]_{(\lambda)}=\langle\gamma,\lambda\rangle$,
$\forall\lambda\in\X^{+}$. If $G$ is simple, then
$\kk[G/\Ru{P}]_0=\kk$, cf.~Lemma~\ref{fixed pt}.

Moreover, the $\gamma$-action defines a vector space grading of
$V(\lambda)$ and a $G$-module grading of $\kk[G]$ such that
$V(\lambda)^{\Ru{P}}$ and $\kk[G/\Ru{P}]_{(\lambda)}$ are the
homogeneous components of maximal degree. In fact, the weight
polytope of $V(\lambda)^{\Ru{P}}$ is the face of the weight
polytope of~$V(\lambda)$, where the linear function
$\langle\gamma,\cdot\rangle$ reaches its maximal value.

It follows that
$\kk[G/\Ru{P}]_{(\lambda)}\cdot\kk[G/\Ru{P}]_{(\mu)}$ is the
homogeneous component of maximal degree in
$\kk[G]_{(\lambda)}\cdot\kk[G]_{(\mu)}$. The isotypic
decomposition of the latter product space is a well-known
particular case of Proposition~\ref{tails}. Taking the maximal
degree means that we must choose only those direct summands with
$\langle\gamma,\lambda+\mu-\beta_i\rangle=
\langle\gamma,\lambda\rangle+\langle\gamma,\mu\rangle$ $\iff$
$\beta_i\perp\gamma$ $\iff$ $\beta_i\in\ZZ_{+}\Pi_L$. The
respective simple $L$-modules are exactly those occurring in the
decomposition of $V(\lambda)^{\Ru{P}}\otimes V(\mu)^{\Ru{P}}$.
Thus the particular case of Proposition~\ref{tails} implies the
general one.
\end{remark}

\begin{remark}
Since $G/\Ru{P}$ is a spherical homogeneous space under $G\times
L$, the powerful theory of spherical varieties \cite{var.sph},
\cite{sph} can be applied to the study of its equivariant
embeddings. For instance, it is easy to deduce from
Proposition~\ref{tails} that the \emph{valuation cone} of
$G/\Ru{P}$ equals $-\dual{C}_L$, and the \emph{colors} are
identified with simple coroots of~$G$. Now it follows from the
general theory that \emph{normal} affine $(G\times
L)$-embeddings $G/\Ru{P}\hookrightarrow X$ are in bijection with
convex polyhedral cones generated by $\dual\Pi$ and finitely
many vectors from $-\dual{C}_L$, $(G\times L)$-orbits in $X$
correspond to faces of such a cone with interiors intersecting
$-\dual{C}_L$, etc.

However, in this paper we prefer to give a more elementary
treatment of affine embeddings of $G/\Ru{P}$ based on properties
of their coordinate algebras and on explicit embeddings into
ambient vector spaces. Our approach is similar to that of
\cite{grp.comp}, \cite[3.3--3.4]{eq.emb} for projective group
completions and reductive monoids.
\end{remark}

\subsection{Affine embeddings}\label{affine.emb}

Affine $(G\times L)$-embeddings $X\hookleftarrow G/\Ru{P}$ are
determined by their coordinate algebras $\kk[X]$, which are
$(G\times L)$-stable finitely generated subalgebras of
$\kk[G/\Ru{P}]$ with the quotient field $\kk(G/\Ru{P})$. Since
$\kk[G/\Ru{P}]$ is multiplicity free, we have
$\kk[X]=\bigoplus_{\lambda\in S}\kk[G/\Ru{P}]_{(\lambda)}$,
where $S\subseteq\X^{+}$ is a finitely generated semigroup such
that $\ZZ{S}=\X(T)$. Proposition~\ref{tails} implies that all
highest weights $\lambda+\mu-\beta$ of $V_L(\lambda)\otimes
V_L(\mu)$ belong to $S$ whenever $\lambda,\mu\in S$.

The variety $X$ is normal iff $\kk[X]^{U\times(U\cap L)}$ is
integrally closed \cite[Thm.\,III.3.3-2]{inv}. But the latter
algebra is just the semigroup algebra of~$S$, which is
integrally closed iff $S=\Sigma^{+}\cap\X(T)$ is the semigroup
of all lattice vectors in the polyhedral cone
$\Sigma^{+}=\QQ_{+}S$. For example, for $X=\CE(G/\Ru{P})$ we
have $S=\X^{+}=C\cap\X(T)$, $\Sigma^{+}=C$.

By Proposition~\ref{tails}, a $(G\times L)$-stable subspace
$\kk[G/\Ru{P}]_{(\lambda_1)}\oplus\dots\oplus\kk[G/\Ru{P}]_{(\lambda_m)}$
generates $\kk[X]$ iff $S$ is $L$-generated by
$\lambda_1,\dots,\lambda_m$ in the sense of the following

\begin{definition}
We say that $\lambda_1,\dots,\lambda_m$ \emph{$L$-generate} $S$
if $S$ consists of all highest weights
$k_1\lambda_1+\dots+k_m\lambda_m-\beta$ of $L$-modules
$V_L(\lambda_1)^{\otimes k_1}\otimes\dots\otimes
V_L(\lambda_m)^{\otimes k_m}$, $k_1,\dots,k_m\in\ZZ_{+}$.
(In particular, any generating set $L$-generates~$S$.)
\end{definition}

Since $\kk[G/\Ru{P}]_{(\lambda)}^{*}\cong
\Hom\bigl(V(\lambda)^{\Ru{P}},V(\lambda)\bigr)$, we have a
$(G\times L)$-equi\-vari\-ant closed immersion
$X\hookrightarrow\bigoplus_{i=1}^m
\Hom\bigl(V(\lambda_i)^{\Ru{P}},V(\lambda_i)\bigr)$ in this case.
Conversely, any equivariant immersion
\begin{equation*}
G/\Ru{P}\hookrightarrow
\bigoplus_{i=1}^m\Hom\bigl(V(\lambda_i)^{\Ru{P}},V(\lambda_i)\bigr)
\end{equation*}
gives rise to an affine embedding
\begin{equation*}
X=\overline{G/\Ru{P}}\subseteq
\bigoplus_{i=1}^m\Hom\bigl(V(\lambda_i)^{\Ru{P}},V(\lambda_i)\bigr)
\end{equation*}
with weight semigroup $S$ $L$-generated by
$\lambda_1,\dots,\lambda_m$.

Put $V=V(\lambda_1)\oplus\dots\oplus V(\lambda_m)$, then
$V^{\Ru{P}}\cong V_L(\lambda_1)\oplus\dots\oplus
V_L(\lambda_m)$. Since a multiple of each $\mu\in C\cap\conv
W_L\{\lambda_1,\dots,\lambda_m\}$ eventually occurs as an
$L$-highest weight in $(V^{\Ru{P}})^{\otimes k}$ by
\cite[Lemma~1]{grp.comp}, we have $\Sigma^{+}=
C\cap\QQ_{+}(W_L\{\lambda_1,\dots,\lambda_m\})$.

We sum up the above discussion in the following theorem:

\begin{theorem}\label{aff.emb}
There is a bijection between affine $(G\times L)$-equivariant
embeddings $X\hookleftarrow G/\Ru{P}$ and subsemigroups
$S\subseteq\X^{+}$ $L$-generated by finitely many weights
$\lambda_1,\dots,\lambda_m\in\X^{+}$ and such that
$\ZZ{S}=\X(T)$. There is a natural equivariant embedding
\begin{align*}
X&\hookrightarrow\bigoplus_{i=1}^m
\Hom\bigl(V(\lambda_i)^{\Ru{P}},V(\lambda_i)\bigr)\\
e\Ru{P}&\mapsto
\bigl(\1_{V(\lambda_1)^{\Ru{P}}},\dots,\1_{V(\lambda_m)^{\Ru{P}}}\bigr)
\end{align*}
for any $L$-generating set of~$S$. The convex cone spanned by
$S$ is
$\Sigma^{+}=C\cap\QQ_{+}(W_L\{\lambda_1,\dots,\lambda_m\})$. The
variety $X$ is normal iff $S=\Sigma^{+}\cap\X(T)$.
\end{theorem}

\begin{example}
If $G$ is semisimple simply connected, then there is a natural
inclusion
\begin{equation*}
\CE(G/\Ru{P})\subseteq\bigoplus_{i=1}^l
\Hom\bigl(V(\omega_i)^{\Ru{P}},V(\omega_i)\bigr)
\end{equation*}
where
$\omega_1,\dots,\omega_l$ are the fundamental weights of~$G$.
\end{example}

\subsection{Relation to reductive monoids}\label{emb<->red.mon}

One observes that the classification of affine embeddings of
$G/\Ru{P}$ is given in the same terms as the classification of
algebraic monoids with the group of invertibles~$L$ \cite{vi},
\cite[3.3]{eq.emb}. Here is a geometric explanation to this
coincidence.

The group $L$ embeds in $G/\Ru{P}$ as the orbit of $e\Ru{P}$.
Let $M$ be the closure of $L$ in~$X$. Under the embedding
$X\hookrightarrow\Hom(V^{\Ru{P}},V)$, $M$~embeds in $\End
V^{\Ru{P}}$ as an algebraic submonoid with the group of
invertibles~$L$. As $V^{\Ru{P}}\cong\bigoplus_{i=1}^m
V_L(\lambda_i)$, we have $\kk[M]=\bigoplus_{\lambda\in S}
\kk[L]_{(\lambda)}$, where $S$ is the semigroup $L$-generated by
$\lambda_1,\dots,\lambda_m$.

There is a natural birational proper map $G\times^PM\to X$,
$(g,z)\mapsto gz$, where $P$ acts on $M$ through its quotient
group $L\cong P/\Ru{P}$ by left translations. Moreover, $X$~is
recovered from $M$ as $\Spec\kk[G\times^PM]$. Indeed,
$\kk[G\times^PM]=(\kk[G]\otimes\kk[M])^P=
\bigoplus_{\lambda,\mu}\bigl(\kk[G]_{(\mu)}\otimes
\kk[M]_{(\lambda)}\bigr)^P$, where $P$ acts on $G$ by right
translations and on $M$ as above. But $\bigl(\kk[G]_{(\mu)}\otimes
\kk[M]_{(\lambda)}\bigr)^P\cong V(\mu)^{*}\otimes\bigl(V(\mu)\otimes
V_L(\lambda)^{*}\bigr)^P\otimes V_L(\lambda)\cong
V(\lambda)^{*}\otimes V_L(\lambda)$ whenever $\lambda=\mu$,
and~$0$, otherwise. Therefore $\kk[G\times^PM]=\kk[X]$.

Conversely, let $M\hookleftarrow L$ be any algebraic monoid, and
$S\subseteq\X^{+}_L$ its weight semigroup. The same reasoning as
above shows that $\kk[G\times^PM]\cong\bigoplus_{\lambda\in
S\cap C}\kk[G/\Ru{P}]_{(\lambda)}$. In other words, affine
embeddings of $G/\Ru{P}$ correspond to algebraic monoids, whose
group of invertible elements is $L$ and the weight semigroup
consists of $G$-dominant weights.

There is a bijective correspondence between the following sets:
\{$(G\times L)$-orbits in~$X$\}, \{$(G\times L)$-stable prime
ideals in $\kk[X]$\}, \{faces $\Gamma\subseteq\Sigma^{+}$ such
that $\nu\notin\Gamma$ whenever $V_L(\nu)\hookrightarrow
V_L(\lambda)\otimes V_L(\mu)$, $\lambda\in S$, $\mu\in
S\setminus\Gamma$\}, \{$(L\times L)$-stable prime ideals in
$\kk[M]$\}, \{$(L\times L)$-orbits in~$M$\}. Thus we have
finally proved

\begin{proposition}\label{emb<->mon}
There is a bijection between affine $(G\times L)$-embeddings
$G/\Ru{P}\hookrightarrow X$ and algebraic monoids $M$ with the
group of invertibles $L$ and the weight semigroup
$S\subseteq\X^{+}_G$, given by $X=\Spec\kk[G\times^PM]$. The
natural proper birational map $G\times^PM\to X$ yields a
bijection between $(G\times L)$-orbits in $X$ and $(L\times
L)$-orbits in $M$ preserving inclusions of orbit closures.
\end{proposition}

\begin{example}
Let $G=\GL(n)$ and $P$ be the stabilizer of a $d$-subspace
in~$\kk^n$. Then $G/\Ru{P}$ embeds in the variety of complexes
\begin{equation*}
X=\left\{\kk^d\stackrel{A_1}{\longrightarrow}
\kk^n\stackrel{A_2}{\longrightarrow}\kk^{n-d}\Bigm|
A_2A_1=0\right\}
\end{equation*}
so that $e\Ru{P}\mapsto(A_1^0,A_2^0)$, where $A_1^0$ is
the inclusion and $A_2^0$ the projector w.r.t.\ a fixed
decomposition $\kk^n=\kk^d\oplus\kk^{n-d}$. Here
$L=\GL(d)\times\GL({n-d})$ is the stabilizer of this
decomposition.

In the above notation, we may take $V=\kk^n\oplus(\kk^n)^{*}$,
so that $V^{\Ru{P}}=\kk^d\oplus(\kk^{n-d})^{*}$, and
$X\hookrightarrow\Hom(V^{\Ru{P}},V)$, $(A_1,A_2)\mapsto
A_1\oplus A_2^{*}$. It follows that $M=\Mat(d)\times\Mat(n-d)=
\{(A_1,A_2)\mid\Im A_1\subseteq\kk^d\subseteq\Ker A_2\}$. The
map $G\times^PM\to X$ is given by
$(g,A_1,A_2)\mapsto(gA_1,A_2g^{-1})$. The weight semigroup $S$
is freely generated by
$\pi_1,\dots,\pi_d,\pi_1^{*},\dots,\pi_{n-d}^{*}$, where $\pi_i$
is the highest weight of $\bigwedge^i\kk^n$ and $\pi_i^{*}$ is
the dual highest weight. The $(G\times L)$-orbits in~$X$, as
well as $(L\times L)$-orbits in~$M$, are determined by the
numerical invariant $(\rk A_1,\rk A_2)$, and the inclusion of
orbit closures corresponds to the product order on these pairs.
\end{example}

\subsection{Orbits}\label{orbits}

Our aim is to describe the orbital decomposition of
$X\hookleftarrow G/\Ru{P}$. Let us recall some basic notions and
introduce some notation.

\begin{definition}
The \emph{generic modality} of the action $G:X$ of an algebraic
group on an irreducible variety is the number
\begin{align*}
d_G(X)&=\codim_XGx\text{ ($x\in X$ a general point)}\\
&=\min_{x\in X}\codim_XGx=\trdeg\kk(X)^G
\end{align*}
The \emph{modality} of $G:X$ is the maximal number of parameters
which a continuous family of $G$-orbits in $X$ depends on, i.e.,
\begin{equation*}
\md_GX=\max_{Y\subseteq X}d_G(Y)
\end{equation*}
where $Y$ runs over all $G$-stable irreducible subvarieties
of~$X$. (Note that $X$ has finitely many orbits iff $\md_GX=0$.)
\end{definition}

Let $\Sigma=\QQ_{+}(W_L\{\lambda_1,\dots,\lambda_m\})$ be the
convex cone generated by the weight polytope of~$V^{\Ru{P}}$.
Note that $\Sigma=W_L\Sigma^{+}$, $\Sigma^{+}=\Sigma\cap C$.

For any face $\Gamma\subseteq\Sigma$, let $V_{\Gamma}\subseteq
V^{\Ru{P}}$ be the sum of $T$-eigenspaces with eigenweights
in~$\Gamma$, and $e_{\Gamma}:V^{\Ru{P}}\to V_{\Gamma}$ be the
$T$-equivariant projector.

For any subset $\Phi\subseteq\Pi$, let $L_{\Phi}$ denote the
standard Levi subgroup with the system of simple roots~$\Phi$,
$L_{\Phi}'$ its commutator subgroup, $P_{\Phi}$ the standard
parabolic subgroup generated by $L_{\Phi}$ and~$B$, and
$P_{\Phi}^{-}$ the opposite parabolic subgroup. If
$N\subseteq\X(T)\otimes\QQ$ is a subspace such that there exists
$\gamma\in\dual{C}$, $\gamma\perp N$, $\gamma\not\perp\alpha$,
$\forall\alpha\in\Delta\setminus N$, then $\Phi=\Pi\cap N$ is
the base of the root subsystem $\Delta\cap N$, and we put
$L_N=L_{\Phi}$, etc. For any sublattice $\Lambda\subseteq\X(T)$,
denote by $T^{\Lambda}\subseteq T$ the diagonalizable group
which is the common kernel of all characters
$\lambda\in\Lambda$.

Suppose that $\Gamma$ is a face of $\Sigma$ whose interior
intersects~$C$. Put
$|\Gamma|=\langle\Delta_L\rangle\cap\langle\Gamma\rangle$,
$\|\Gamma\|=|\Gamma|\oplus\langle\Gamma\rangle^{\perp}$,
$\langle\Gamma\rangle_{\ZZ}=\sum_{\lambda_i\in\Gamma}\ZZ\lambda_i
+(\ZZ\Delta_L\cap\langle\Gamma\rangle)$ (a sublattice
generating~$\langle\Gamma\rangle$).

The following theorem is a counterpart of the results of
\cite[\S9]{grp.comp}.

\begin{theorem}\label{orb}
The $(G\times L)$-orbits $Y\subseteq X$ are in bijection with
the faces $\Gamma\subseteq\Sigma$ whose interiors intersect~$C$.
The inclusion of faces corresponds to the inclusion of orbit
closures. The orbit $Y=Y_{\Gamma}$ is represented
by~$e_{\Gamma}$. The stabilizers are:
\begin{align*}
(G\times L)_{e_{\Gamma}}&=
\left[\Ru{(P_{\|\Gamma\|})}\times\Ru{\bigl(L\cap
P_{\|\Gamma\|}^{-}\bigr)}\right]\\
&\qquad\leftthreetimes
\left[\Bigl(L_{\langle\Gamma\rangle^{\perp}}'
T^{\langle\Gamma\rangle_{\ZZ}}
\times\bigl(L\cap L_{\langle\Gamma\rangle^{\perp}}\bigr)'
T^{\langle\Gamma\rangle_{\ZZ}}\Bigr)
\cdot\diag L_{|\Gamma|}\right]\\
G_{e_{\Gamma}}&=\Ru{(P_{\|\Gamma\|})}\leftthreetimes
L_{\langle\Gamma\rangle^{\perp}}'T^{\langle\Gamma\rangle_{\ZZ}}
\end{align*}
All $G$-orbits in $Y$ are isomorphic and permuted transitively
by~$L$. The (generic) modality of $Y$ is:
\begin{equation*}
d_G(Y)=\dim L\cap\Ru{(P_{\|\Gamma\|})}
\end{equation*}
\end{theorem}

\begin{proof}
By Proposition~\ref{emb<->mon}, $(G\times L)$-orbits in $X$ are
in bijection with $(L\times L)$-orbits in
$M=\overline{L}\subseteq\End V^{\Ru{P}}$. Therefore it suffices
to describe the orbits for $L\times L:M$. This description goes
back to Putcha and Renner. In particular, one finds out that the
projectors $e_{\Gamma}$ form a complete set of orbit
representatives (cf.~\cite[Thm.\,8]{grp.comp}). Let us give an
outline of an elementary proof.

First observe that $\overline{T}$ intersects all $(L\times
L)$-orbits in~$M$. (For $\kk=\mathbb{C}$, the easiest way to see
it is to close in $M$ the Cartan decomposition $L=KTK$, where
$K\subset L$ is a maximal compact subgroup. For arbitrary~$\kk$,
one may consider the Iwahori decomposition of
$G\bigl(\kk((t))\bigl)$ instead, see \cite[2.4, Exemple~2]{var.sph}.)
Next, it is easy to deduce from affine toric geometry that
$T$-orbits in $\overline{T}$ are represented by $e_{\Gamma}$
over \emph{all} faces $\Gamma\subseteq\Sigma$.

But one sees from the structure of $(L\times L)_{e_{\Gamma}}$
that $\bigl((L\times L)e_{\Gamma}\bigr)^{\diag{T}}$ is a union
of $T$-orbits permuted by $W_L$ transitively. Indeed, if
$y=(g_1,g_2)e_{\Gamma}$ is fixed by $\diag{T}$, then one may
assume that $(g_1,g_2)^{-1}(\diag{T})(g_1,g_2)$ is contained in
the maximal torus $\bigl(T^{\langle\Gamma\rangle_{\ZZ}}\times
T^{\langle\Gamma\rangle_{\ZZ}}\bigr)\cdot\diag{T}$ of $(L\times
L)_{e_{\Gamma}}$. Hence $g_1,g_2\in N_L(T)$ represent two
elements $w_1,w_2\in W_L$ acting on $\Gamma$ equally, and
$y=w_2e_{\Gamma}w_2^{-1}$. Thus $(L\times L)$-orbits in~$M$ are
represented by those $e_{\Gamma}$ corresponding to faces with
interiors intersecting~$C$.

The above reasoning also proves the assertion on inclusions of
faces and orbit closures, since it is true for $T$-orbits
in~$\overline{T}$.

Now we compute the stabilizers. Let $V_{\Gamma}'$ be the
$T$-stable complement to $V_{\Gamma}$ in~$V^{\Ru{P}}$. For
$(g,h)\in G\times L$ we have: $ge_{\Gamma}h^{-1}=e_{\Gamma}$
iff
\begin{enumerate}
\renewcommand{\theenumi}{\textup{(\arabic{enumi})}}
\renewcommand{\labelenumi}{\theenumi}
\item\label{Im} $gV_{\Gamma}=V_{\Gamma}$,
\item\label{Ker} $hV_{\Gamma}'=V_{\Gamma}'$,
\item\label{diag} the actions of $g,h$ on $V_{\Gamma}\cong
V^{\Ru{P}}/V_{\Gamma}'$ coincide.
\end{enumerate}

The condition~\ref{Im} means that $g\in P_{\|\Gamma\|}$. Indeed,
for any $\alpha\in\Delta^{+}\setminus\|\Gamma\|$ we have
$\langle\dual\alpha,\Gamma\rangle\ge0$, and the strict
inequality is achieved. Hence $\Ru{(P_{\|\Gamma\|})}$ fixes
$V_{\Gamma}$ pointwise, whereas no element of
$\Ru{\bigl(P_{\|\Gamma\|}^{-}\bigr)}$ preserves~$V_{\Gamma}$. By
definition, $V_{\Gamma}$ is $L_{|\Gamma|}$-stable.
On the other hand, it is easy to see that
adding roots $\alpha\in\|\Gamma\|\setminus|\Gamma|$ moves the
weights of $V(\lambda_i)\cap V_{\Gamma}$ outside the weight
polytope of~$V(\lambda_i)$, $\forall i=1,\dots,m$. Hence the
respective root vectors act on $V_{\Gamma}$ trivially, i.e.,
$\alpha\perp\Gamma$. This means that
$\Delta\cap\|\Gamma\|=(\Delta\cap|\Gamma|)\sqcup
(\Delta\cap\langle\Gamma\rangle^{\perp})$ is a disjoint
orthogonal union, $L_{\|\Gamma\|}=L_{|\Gamma|}\cdot
L_{\langle\Gamma\rangle^{\perp}}$, and
$L_{\langle\Gamma\rangle^{\perp}}'$ fixes $V_{\Gamma}$ pointwise.

Similar arguments show that \ref{Ker}~$\iff h\in
L\cap P_{\|\Gamma\|}^{-}$, and the subgroup $\Ru{\bigl(L\cap
P_{\|\Gamma\|}^{-}\bigr)}\leftthreetimes
\bigl(L\cap L_{\langle\Gamma\rangle^{\perp}}\bigr)'$ acts on
$V^{\Ru{P}}/V_{\Gamma}'$ trivially. Thus after factoring out the
kernels of the actions, we may assume $g,h\in
L_{|\Gamma|}\subseteq L$. But $L_{|\Gamma|}$ acts on
$V_{\Gamma}\cong V^{\Ru{P}}/V_{\Gamma}'$ with kernel
$T^{\langle\Gamma\rangle_{\ZZ}}$,
$\langle\Gamma\rangle_{\ZZ}$~being the weight lattice
of~$V_{\Gamma}$. Hence \ref{diag} $\iff g\equiv h\mod
T^{\langle\Gamma\rangle_{\ZZ}}$, and we are done.

The formula for $G_{e_{\Gamma}}$ stems from that for
$(G\times L)_{e_{\Gamma}}$ immediately. Since the $L$-action on
$Y$ commutes with the $G$-action, it permutes the $G$-orbits
transitively, and all of them are isomorphic and, in particular,
have the same dimension. Now $d_G(Y)=\dim Y-\dim
Ge_{\Gamma}=\dim(G\times L)-\dim(G\times L)_{e_{\Gamma}}-\dim G
+\dim G_{e_{\Gamma}}=\dim L-
\dim\Ru{\bigl(L\cap P_{\|\Gamma\|}^{-}\bigr)}-
\dim\bigl(L\cap L_{\langle\Gamma\rangle^{\perp}}\bigr)'
T^{\langle\Gamma\rangle_{\ZZ}}-
\dim L_{|\Gamma|}/T^{\langle\Gamma\rangle_{\ZZ}}=
\dim L\cap\Ru{(P_{\|\Gamma\|})}$.
\end{proof}

\begin{corollary}\label{orb.can}
The $(G\times L)$-orbits $Y\subseteq\CE(G/\Ru{P})$ are in
bijection with the subsystems of simple roots
$\Pi_Y\subseteq\Pi$ such that no component of $\Pi_Y$ is
contained in~$\Pi_L$. The stabilizers of $Y$ in $G\times L$ and
in $G$ are:
\begin{gather*}
\left[\Ru{\bigl(P_{\Pi_Y\cup(\Pi_L\setminus\partial\Pi_Y)}\bigr)}
\times\Ru{\bigl(L\cap
P_{\Pi_L\setminus\partial\Pi_Y}^{-}\bigr)}\right]\qquad\qquad\\
\qquad\qquad\leftthreetimes
\left[\Bigl(L_{\Pi_Y}'
\times\bigl(L\cap L_{\Pi_Y}'\bigr)\Bigr)
\cdot\diag L_{\Pi_L\cap\Pi_Y^{\perp}}\right]\\
\text{and\qquad}
\Ru{\bigl(P_{\Pi_Y\cup(\Pi_L\setminus\partial\Pi_Y)}\bigr)}
\leftthreetimes L_{\Pi_Y}'
\end{gather*}
Here $\partial\Pi_Y$ is the set of simple roots from
$\Pi\setminus\Pi_Y$ neighboring with $\Pi_Y$ on the Dynkin
diagram of~$G$. We have $d_G(Y)=\dim
L\cap\Ru{\bigl(P_{\Pi_L\setminus\partial\Pi_Y}\bigr)}$.

The modality $\md_G\CE(G/\Ru{P})=\max_Yd_G(Y)$ is reached on $Y$ such that
$\Pi_Y\supseteq\Pi\setminus\Pi_L$, and each component of $\Pi_Y$
is obtained from a component of $\Pi\setminus\Pi_L$ by adding
roots from $\Pi_L$ in such a way that $\dim
L_{\Pi_L^{(k)}\setminus\partial\Pi_Y}=\min$ for each component
$\Pi_L^{(k)}\subseteq\Pi_L$. (In particular,
$\md_G\CE(G/\Ru{P})=0$ iff $\Pi_L$ is a union of components
of~$\Pi$, which implies Proposition~\ref{pr2}.)
\end{corollary}

\begin{proof}
The orbits $Y\subseteq\CE(G/\Ru{P})$ are in bijection with the
faces $\Gamma\subseteq W_LC$ whose interiors intersect~$C$. Then
$\Gamma\cap C$ is the face of $C$ of the same dimension
$\dim\Gamma$, and the dual face of $\QQ_{+}\dual\Pi$ is spanned
by a certain subset $\dual\Pi_Y\subseteq\dual\Pi$, so that
$\langle\beta_Y,W_LC\rangle\ge0$ for a certain positive
linear combination $\beta_Y$ of~$\dual\Pi_Y$.

Suppose that a component $\Pi_Y^{(i)}$ of $\Pi_Y$ is contained
in~$\Pi_L$. Let $\rho_{(i)}\in C$ be the sum of fundamental
weights corresponding to simple roots in~$\Pi_Y^{(i)}$, and
$w_{(i)}\in W_L$ the longest element of the Weyl group
of~$\Pi_Y^{(i)}$. Then $w_{(i)}\rho_{(i)}\in W_LC$, but
$\langle\beta_Y,w_{(i)}\rho_{(i)}\rangle=
-\langle\beta_Y,\rho_{(i)}\rangle<0$, a contradiction.

On the other hand, if no component of $\Pi_Y$ is contained
in~$\Pi_L$, then it is easy to find a positive linear
combination $\beta_Y$ of~$\dual\Pi_Y$ such that
$\langle\beta_Y,\Pi_L\rangle\le0$ (choosing sufficiently large
coefficients in $\beta_Y$ successively at
$\dual\alpha\in\dual\Pi_Y\cap\dual\Pi_L$ and finally at
$\dual\alpha\in\dual\Pi_Y\setminus\dual\Pi_L$), whence
$\langle\beta_Y,W_LC\rangle\ge0$.

We have $|\Gamma|=\langle\Pi_L\rangle\cap\Pi_Y^{\perp}$,
$\langle\Gamma\rangle^{\perp}=\langle\Pi_Y\rangle$, whence
$\|\Gamma\|$ is spanned by $\Pi_Y\sqcup(\Pi_L\cap\Pi_Y^{\perp})=
\Pi_Y\cup(\Pi_L\setminus\partial\Pi_Y)$. Choosing a sufficiently
large generating set $\{\lambda_1,\dots,\lambda_m\}$
for~$\X^{+}$, we see that
$\langle\Gamma\rangle_{\ZZ}=\langle\Gamma\rangle\cap\X(T)$,
whence $T^{\langle\Gamma\rangle_{\ZZ}}$ is connected and is in
fact a maximal torus in~$L_{\Pi_Y}'$. This proves the formul{\ae}
for the stabilizers and for~$d_G(Y)$.

Finally, if we look for an $Y$ with $d_G(Y)=\max$, i.e.,
$\dim L_{\Pi_L\setminus\partial\Pi_Y}=\min$, we may always
include $\Pi\setminus\Pi_L$ in $\Pi_Y$ in order to enlarge
$\Pi_L\cap\partial\Pi_Y$ as much as possible. It remains to note
that $\dim L_{\Pi_L\setminus\partial\Pi_Y}=\min\iff
\dim L_{\Pi_L^{(k)}\setminus\partial\Pi_Y}=\min$, $\forall k$.
\end{proof}

\begin{example}\label{Mat(n,n-1)}
Let $G=\SL(n)$ and $P$ be the stabilizer of a hyperplane
in~$\kk^n$. Then $\CE(G/\Ru{P})=\Mat(n,n-1)$ with the $G$-action
by left multiplication. Let $\alpha_1,\dots,\alpha_{n-1}$ be the
simple roots of~$G$. The group $L\cong\GL(n-1)$ acts on
$\Mat(n,n-1)$ by right multiplication,
$\Pi_L=\{\alpha_1,\dots,\alpha_{n-2}\}$.

The possible choices for $\Pi_Y$ are
$\Pi_Y=\{\alpha_k,\dots,\alpha_{n-1}\}$, $1\le k\le n$. The
respective orbit $Y$ consists of all matrices of rank $k-1$. We
have $L_{\Pi_L\setminus\partial\Pi_Y}\cong\GL(k-1)\times
\GL(n-k)$, and $d_G(Y)=(k-1)(n-k)$. The latter formula can be
derived directly from the observation that a $G$-orbit in $Y$ is
formed by all matrices of rank $k-1$ with given linear
dependencies between the columns. The space of linear
dependencies depends on $(k-1)(n-k)$ parameters, which are the
coefficients of linear expressions of all columns through the
basic ones. (In fact, the respective moduli space is nothing
else, but the Grassmannian of ${(n-k)}$-subspaces in
$\kk^{n-1}$.) The maximal value of $d_G(Y)$ is reached for
$k=[(n+1)/2]$, and $\md_{\SL(n)}\Mat(n,n-1)=[(n-1)^2/4]$.
\end{example}

We illustrate all the above results by another example of a
canonical embedding:

\begin{example}\label{iso.Im}
Let $G=\Sp(2l)$ and $P$ be the stabilizer of a Lagrangian
subspace $\kk^l\subset\kk^{2l}$. A complementary Lagrangian
subspace is canonically isomorphic to $(\kk^l)^{*}$, the pairing
with $\kk^l$ being given by the symplectic form. Then $L=\GL(l)$
is the stabilizer of the decomposition
$\kk^{2l}=\kk^l\oplus(\kk^l)^{*}$.

Let $X\subseteq\Mat(2l,l)$ be the set of all linear maps with
isotropic image. Then $X$ is an affine embedding of~$G/\Ru{P}$,
so that $e\Ru{P}$ is mapped to the identity map
$\kk^l\to\kk^l\subset\kk^{2l}$. In the notation
of~\ref{emb<->red.mon}, we have $M=\Mat(l)$, and $G\times^PM\to
X$, $(g,A)\mapsto g\genfrac(){}{}A0$, is the
multiplication map.

Let $\alpha_1,\dots,\alpha_l$ be the simple roots and
$\omega_1,\dots,\omega_l$ the fundamental weights of~$G$, in the
standard order. Then $\Pi_L=\{\alpha_1,\dots,\alpha_{l-1}\}$,
and~$\omega_i$, considered as a dominant weight of~$L$, is the
highest weight of~$\bigwedge^i\kk^l$, $\forall i$. It follows
that $S$ is generated by $\omega_1,\dots,\omega_l$, since
$S\ni\omega_1$. Therefore $X=\CE(G/\Ru{P})$.

The $(G\times L)$-orbits $Y\subset X$ (and $(L\times L)$-orbits
in~$M$) consist of all matrices of given rank~$k$, $0\le k\le
l$. We have $\Pi_Y=\{\alpha_{k+1},\dots,\alpha_l\}$,
$L_{\Pi_L\setminus\partial\Pi_Y}\cong\GL(k)\times \GL(l-k)$,
$d_G(Y)=k(l-k)$, and $\md_GX=[l^2/4]$. The reasoning is similar
to that of Example~\ref{Mat(n,n-1)}.
\end{example}

\subsection{Smoothness}\label{smooth.emb}

Now we classify those affine embeddings of $G/\Ru{P}$ which are
smooth.

\begin{example}\label{3:smooth}
Here are three basic examples of smooth embeddings
$X\hookleftarrow G/\Ru{P}$:
\begin{enumerate}
\renewcommand{\theenumi}{\textup{(\arabic{enumi})}}
\renewcommand{\labelenumi}{\theenumi}
\item\label{SL/P_u->Mat} The embedding $X=\Mat(n,n-1)$ of
Example~\ref{Mat(n,n-1)} is smooth.
\item\label{GL->Mat} The embedding $X=\Mat(n)$ of $G=\GL(n)$ is
smooth. (Here $P=L=G$.)
\item\label{X=G} The group $G$ itself is a smooth embedding
of~$G$. (Again $P=L=G$.)
\end{enumerate}
\end{example}

Our next result shows that these are the only nontrivial
examples of smooth affine embeddings.

\begin{theorem}\label{smooth}
Any smooth affine $(G\times L)$-embedding of $G/\Ru{P}$ is of
the form $X=G_0\times^{Z_0}X_{\perp}$. Here $G=(G_0\times
G_{\perp})/Z_0$ is the quotient of a product of two reductive
groups by a finite central diagonally embedded subgroup,
$P\supseteq G_0$, and the embedding
$G_{\perp}/\Ru{P}\hookrightarrow X_{\perp}$ is the direct
product of several embeddings \ref{SL/P_u->Mat}, \ref{GL->Mat}
of Example~\ref{3:smooth}, where the actions on the factors of
type~\ref{SL/P_u->Mat} are possibly shifted by some characters
of~$G_{\perp}$.
\end{theorem}

\begin{proof}
The idea of the proof is similar to that of
\cite[Thm.\,9]{grp.comp}.

Let $\Gamma_0$ be the minimal face (i.e., maximal linear
subspace) in~$\Sigma$. Then $\langle\gamma,\Gamma_0\cap
C\rangle\ge0$ for $\gamma$ as in Remark~\ref{grading}. As
$\gamma$ is fixed by~$W_L$, we obtain
$\langle\gamma,\Gamma_0\rangle\ge0\implies\gamma\perp\Gamma_0$.
It follows that $\Gamma_0$ is orthogonal to each component of
$\Pi$ not contained in~$\Pi_L$.

Therefore each component of $\Pi$ is either contained in
$|\Gamma_0|$ or orthogonal to~$\Gamma_0$. Since $X$ is smooth,
hence normal, $\langle\Gamma_0\rangle_{\ZZ}=\Gamma_0\cap\X(T)$,
whence $T_{\perp}=T^{\langle\Gamma_0\rangle_{\ZZ}}$ is a torus.
Put $G_{\perp}=L_{\langle\Gamma_0\rangle^{\perp}}'T_{\perp}$,
and $G_0=L_{|\Gamma_0|}'T_0$, where the subtorus $T_0\subseteq
T$ extends the maximal torus of $L_{|\Gamma_0|}'$ in such a way
that $T=(T_0\times T_{\perp})/\diag{Z_0}$, $Z_0=T_0\cap
T_{\perp}$ being a finite central subgroup of~$G$.
Then $G=(G_0\times G_{\perp})/\diag{Z_0}$.

It may happen that $\Sigma\cap\gamma^{\perp}\ne\Gamma_0$.
However, the interior of the cone dual to
$\Sigma\cap\gamma^{\perp}$ contains a nonzero vector
$\gamma_0\perp\Pi$: otherwise this interior is separated from
$\Pi^{\perp}$ by a linear function $\langle\cdot,\lambda\rangle$
for some $\lambda\in(\Sigma\cap\gamma^{\perp}\setminus\Gamma_0)
\cap\langle\Pi\rangle\subseteq
\Sigma\cap\langle\Pi_L\rangle\setminus\Gamma_0=\emptyset$, a
contradiction. Replacing $\gamma$ by a multiple of
$\gamma+\gamma_0$, $\gamma_0$~sufficiently small, we may assume
that $\gamma\in\X_{*}(Z)$ and
$\langle\gamma,\Sigma\setminus\Gamma_0\rangle>0$.

By Theorem~\ref{orb}, the face $\Gamma_0$ corresponds to the
closed $(G\times L)$-orbit $Y_0\ni e_{\Gamma_0}$.
The $\gamma$-action by right translations of an argument yields
an equivariant retraction $X\to Y_0$,
$x\mapsto\lim_{t\to\infty}\gamma(t)*x$,
$e_{\Sigma}\mapsto e_{\Gamma_0}$. Thus
we have
$X=(G\times L)\times^{(G\times L)_{e_{\Gamma_0}}}X_{\perp}$,
$X_{\perp}=\overline{(G\times L)_{e_{\Gamma_0}}e_{\Sigma}}\ni
e_{\Gamma_0}$. But $(G\times L)_{e_{\Gamma_0}}=
\bigl(G_{\perp}\times(L\cap G_{\perp})\bigr)\cdot\diag{G_0}$,
and $\diag{G_0}$ acts on $X_{\perp}$ trivially as a normal
subgroup in the stabilizer of the open orbit. Thus
$X=G_0\times^{Z_0}X_{\perp}$, where $X_{\perp}\hookleftarrow
G_{\perp}/\Ru{P}$ is an embedding with a fixed point.

In the sequel, we may assume that $X\subseteq
\bigoplus_{i=1}^m\Hom\bigl(V(\lambda_i)^{\Ru{P}},V(\lambda_i)\bigr)$
itself contains the fixed point~$0$.
The $\gamma$-action contracts the ambient vector space on the
r.h.s.\ to~$0$.
After renumbering the
$\lambda_i$'s, we may assume that
$T_0X=\bigoplus_{i=1}^p\Hom\bigl(V(\lambda_i)^{\Ru{P}},V(\lambda_i)\bigr)$,
$p\le m$. Since $X$ is smooth
and contracted to $0$ by~$\gamma$,
it projects onto $T_0X$ isomorphically.

Let $e\Ru{P}\mapsto(e_1,\dots,e_p)$ under this isomorphism. Then
$e_i$ has the dense $(G\times L)$-orbit in
$\Hom\bigl(V(\lambda_i)^{\Ru{P}},V(\lambda_i)\bigr)$ and commutes
with~$L$, whence by Schur's lemma $e_i$ is a nonzero scalar
operator on $V(\lambda_i)^{\Ru{P}}$. After rescaling the above
isomorphism, we may assume $e_i=\1_{V(\lambda_i)^{\Ru{P}}}$.

Let $G_i\subseteq\GL(V(\lambda_i))$ be the image of~$G$, and
$P_i,L_i$ the images of $P,L$. Then $G_ie_i$ is dense in
$\Hom\bigl(V(\lambda_i)^{\Ru{P}},V(\lambda_i)\bigr)$. It follows that the
orbit of the highest weight vector is dense in~$V(\lambda_i)$,
whence $G_i$ acts on $\PP(V(\lambda_i))$ transitively. By
\cite{Aut(C-flag)}, \cite{Aut(flag)}, $G_i=\GL(V(\lambda_i))$,
$\SL(V(\lambda_i))$, $\Sp(V(\lambda_i))$, or
$\Sp(V(\lambda_i))\cdot\kk^{\times}$, $P_i\ne G_i$ in the 2-nd
case, and $P_i$ fixes the highest weight line in the last two
cases (so that $V(\lambda_i)^{\Ru{P}}\ne V(\lambda_i)$ and $\dim
V(\lambda_i)^{\Ru{P}}=1$, respectively).

Two simple components of $G$ never project to one and the same
$G_i$ non-trivially (because their images must commute).
However, there might exist a simple component of $G$ projecting
to several $G_i$'s non-trivially. Let $i=i_1,\dots,i_q$ be the
respective indices, and
$G_{i_1,\dots,i_q}$, $P_{i_1,\dots,i_q}$, $L_{i_1,\dots,i_q}$ the
images of $G,P,L$ in $G_{i_1}\times\dots\times G_{i_q}$. Then
$G_{i_1,\dots,i_q}'$ is simple, $\dim Z(G_{i_1,\dots,i_q})\le
q$, the orbit $G_{i_1,\dots,i_q}(e_{i_1},\dots,e_{i_q})$ is
dense in
$\bigoplus_{k=1}^q\Hom\bigl(V(\lambda_{i_k})^{\Ru{P}},V(\lambda_{i_k})\bigr)$,
and the stabilizer of $(e_{i_1},\dots,e_{i_q})$ in
$G_{i_1,\dots,i_q}'\cap L_{i_1,\dots,i_q}$ is trivial. In
particular, we have an inequality
\begin{displaymath}
\dim G_{i_1,\dots,i_q}-
\dim\Ru{(P_{i_1,\dots,i_q})}\ge\sum_{k=1}^q\dim V(\lambda_{i_k})
\cdot\dim V(\lambda_{i_k})^{\Ru{P}}
\end{displaymath}
which is strict whenever
$\dim(L_{i_1,\dots,i_q})_{(e_{i_1},\dots,e_{i_q})}>0$. This
leaves the following possibilities:
\begin{enumerate}
\renewcommand{\theenumi}{\textup{(\arabic{enumi})}}
\renewcommand{\labelenumi}{\theenumi}

\item\label{G=P} $G_{i_k}=\GL(n)=P_{i_k}$, $V(\lambda_{i_k})=\kk^n$,
$q+(n^2-1)\ge qn^2$;

\item\label{G>P} $G_{i_k}=\GL(n)\text{ or }\SL(n)$, $V(\lambda_{i_k})=\kk^n$,
$P_{i_k}$~is the stabilizer of the hyperplane in $\kk^n$ given
by vanishing of the last coordinate, $q+(n^2-1)-(n-1)\ge
qn(n-1)$;

\item $G_{i_k}=\GL(n)\text{ or }\SL(n)$,
$V(\lambda_{i_k})=\kk^n\text{ or }(\kk^n)^{*}$ (both cases
occur), $P_{i_k}$~is the stabilizer of the subspace in $\kk^n$
generated by the first $d$ basic vectors, $q+(n^2-1)-d(n-d)\ge
nd+n(n-d)+(q-2)n$.

\end{enumerate}
In all cases we have either $q=1$ or, in the last two cases,
$q=n=2$, and the inequalities become equalities. But in the
latter situation $\dim Z(G_{i_1,\dots,i_q})=q$, and it is easy
to see that
$\dim(L_{i_1,\dots,i_q})_{(e_{i_1},\dots,e_{i_q})}>0$, a
contradiction.

Thus $G\hookrightarrow G_1\times\dots\times G_p$,
$G'=G_1'\times\dots\times G_p'$, and $G(e_1,\dots,e_p)$ is dense
in $\bigoplus_{i=1}^p\Hom\bigl(V(\lambda_i)^{\Ru{P}},V(\lambda_i)\bigr)$,
with stabilizer~$\Ru{P}$. Now an easy dimension count shows that
each triple $(G_i,P_i,V(\lambda_i))$ belongs to case~\ref{G=P}
or \ref{G>P}, and $\dim Z(G)$ is the number of occurrences
of~\ref{G=P}. Thus $X=X_1\times\dots\times X_p$, and each
$X_i=\Hom\bigl(V(\lambda_i)^{\Ru{P}},V(\lambda_i)\bigr)$ is an
embedding of $G_i/\Ru{(P_i)}$, in case~\ref{G=P}, or
$G_i'/\Ru{(P_i)}$, in case~\ref{G>P}.
\end{proof}

\begin{corollary}\label{smooth(can)}
The canonical embedding $\CE(G/\Ru{P})$ is smooth iff
$G=(G_0\times G_{\perp})/Z_0$ is the quotient of a product of
two reductive groups by a finite central diagonally embedded
subgroup, $G_{\perp}=G_1\times\dots\times G_p$, $G_i=\SL(n_i)$
($i>0$), $P\supseteq G_0$, and $P\cap G_i$ are the stabilizers of
hyperplanes (or lines) in~$\kk^{n_i}$.
\end{corollary}

\subsection{Tangent spaces}\label{tan.spaces}

Finally, we shall describe the tangent space $T_0\CE(G/\Ru{P})$
of $\CE(G/\Ru{P})$ at the unique $G$-fixed point~$0$, assuming
that $G$ is simple and $P\ne G$ (see Lemma~\ref{fixed pt}). The
$G$-module structure of this tangent space provides information
on ambient $G$-modules for $\CE(G/\Ru{P})$; namely
$T_0\CE(G/\Ru{P})$ is the smallest one.

As $\kk[\CE(G/\Ru{P})]$ is non-negatively graded by a
one-parameter subgroup $\gamma\in\X_{*}(Z)$ so that
$\kk[\CE(G/\Ru{P})]_0=\kk$ (see Remark~\ref{grading}), the space
$T_0\CE(G/\Ru{P})$ is dual to the linear span of a minimal
system of homogeneous generators for $\kk[\CE(G/\Ru{P})]$. Thus
to describe $T_0\CE(G/\Ru{P})$ is the same thing as to find the
minimal homogeneous generating subspace for
$\kk[\CE(G/\Ru{P})]$, or to find the minimal $L$-generating set
for~$\X^{+}$.

For simplicity, we assume that $G$ is simply connected. Then
$\X^{+}$ is freely generated by the fundamental weights
$\omega_1,\dots,\omega_l$, and it suffices to find out which
$\omega_i$ are $L$-generated by the other fundamental weights.

Let $\alpha_1,\dots,\alpha_l$ be the simple roots of~$G$, and
$\dual\alpha_i$, $\dual\omega_i$ denote the simple coroots and
the fundamental coweights, respectively.

\begin{definition}
The \emph{singularity} of a Dynkin diagram is either the node of
branching or the node representing the long root neighboring
with a short one.
\end{definition}

The $Z$-action by right translations of an argument defines an
invariant algebra multi-grading of $\kk[G/\Ru{P}]$ so that
$\kk[G/\Ru{P}]_{(\lambda)}$ has the weight~$\lambda|_Z$. A
choice of $\gamma\in\X_{*}(Z)\otimes\QQ=\Pi_L^{\perp}=
\langle\dual\omega_i\mid\alpha_i\notin\Pi_L\rangle$ yields a
specialization of this multi-grading, so that
$\deg\kk[G/\Ru{P}]_{(\lambda)}=\langle\gamma,\lambda\rangle$,
cf.~Remark~\ref{grading}. (The degrees might be rational
numbers, however, multiplying $\gamma$ by a sufficiently large
number yields an integer grading.) For brevity, we shall speak
about the \emph{degree} of $\lambda$ w.r.t.~$\gamma$.

Put $\res\lambda=\lambda|_{T\cap L'}$,
$\forall\lambda\in\X^{+}$. Then $\resomega_i$ is a fundamental
weight of the commutator group $L'$ whenever $\alpha_i\in\Pi_L$,
or zero, otherwise. Note that $V_L(\lambda)\hookrightarrow
V_L(\lambda_1)\otimes\dots\otimes V_L(\lambda_n)$ iff
$V(\res\lambda)\hookrightarrow
V(\res\lambda_1)\otimes\dots\otimes V(\res\lambda_n)$ and
$\deg\lambda=\deg\lambda_1+\dots+\deg\lambda_n$
w.r.t.~$\forall\gamma\in\X_{*}(Z)\otimes\QQ$.

The degrees w.r.t.\ the generators $\dual\omega_i$ are
determined in terms of the matrix
$(\langle\dual\omega_i,\omega_j\rangle)_{i,j=1}^l$, which is the
inverse transpose of the Cartan matrix of~$G$. These matrices
are computed in \cite[Table~2]{sem}. The $i$-th row of this
matrix represents the degrees $d_1,\dots,d_l$ of
$\omega_1,\dots,\omega_l$ w.r.t.~$\dual\omega_i$.
Let us label the nodes of the Dynkin diagram by these degrees.
An inspection of the inverse transposed Cartan matrices yields
the following observation:

\sloppy

\begin{itemize}

\item The labels of nodes in a segment from an extreme node up
to either $\alpha_i$ or the singularity form a sequence
$a,2a,\dots,pa$.
\begin{center}
\unitlength 0.40ex
\linethickness{0.4pt}
\begin{picture}(42.00,7.00)
\put(40.00,2.00){\makebox(0,0)[cc]{${\cdots}$}}
\put(31.00,2.00){\line(1,0){4.00}}
\put(30.00,2.00){\circle*{2.00}}
\put(20.00,2.00){\makebox(0,0)[cc]{${\cdots}$}}
\put(2.00,2.00){\circle{2.00}}
\put(29.00,2.00){\line(-1,0){4.00}}
\put(10.00,2.00){\circle{2.00}}
\put(11.00,2.00){\line(1,0){4.00}}
\put(9.00,2.00){\line(-1,0){6.00}}
\put(2.00,4.00){\makebox(0,0)[cb]{$a$}}
\put(10.00,4.00){\makebox(0,0)[cb]{$2a$}}
\put(30.00,3.00){\makebox(0,0)[cb]{$pa$}}
\end{picture}
\end{center}

\item If the Dynkin diagram has no branching, then the nodes
after the singularity up to $\alpha_i$ are labeled by
$da,\dots,da$ or ${(p+1)a/d},a,\dots,a$, where $d$ is the
multiplicity of the ``thick'' edge, depending on whether
$\alpha_i$ is a long root or not.
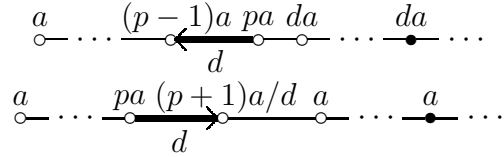
\begin{center}
\unitlength 0.40ex
\linethickness{0.4pt}
\begin{picture}(82.00,10.00)
\put(2.00,4.00){\circle{2.00}}
\put(12.00,4.00){\makebox(0,0)[cc]{${\cdots}$}}
\put(3.00,4.00){\line(1,0){4.00}}
\put(60.00,4.00){\makebox(0,0)[cc]{${\cdots}$}}
\put(42.00,4.00){\circle{2.00}}
\put(51.00,4.00){\line(1,0){4.00}}
\multiput(26.83,4.00)(0.12,0.12){17}{\line(0,1){0.12}}
\multiput(26.83,4.00)(0.12,-0.12){17}{\line(0,-1){0.12}}
\put(69.00,4.00){\line(-1,0){4.00}}
\put(70.00,4.00){\circle*{2.00}}
\put(80.00,4.00){\makebox(0,0)[cc]{${\cdots}$}}
\put(71.00,4.00){\line(1,0){4.00}}
\put(26.00,4.00){\circle{2.00}}
\put(50.00,4.00){\circle{2.00}}
\put(49.00,4.00){\line(-1,0){6.00}}
\put(2.00,6.50){\makebox(0,0)[cb]{$a$}}
\put(27.00,5.00){\makebox(0,0)[cb]{$(p-1)a$}}
\put(42.00,5.50){\makebox(0,0)[cb]{$pa$}}
\put(50.00,6.50){\makebox(0,0)[cb]{$da$}}
\put(70.00,6.50){\makebox(0,0)[cb]{$da$}}
\put(34.11,2.00){\makebox(0,0)[ct]{$d$}}
\put(25.00,4.00){\line(-1,0){8.00}}
\linethickness{2pt}
\put(41.00,4.00){\line(-1,0){14.00}}
\end{picture}
\hfill
\unitlength 0.40ex
\linethickness{0.4pt}
\begin{picture}(89.00,14.00)
\put(2.00,4.00){\circle{2.00}}
\put(12.00,4.00){\makebox(0,0)[cc]{${\cdots}$}}
\put(3.00,4.00){\line(1,0){4.00}}
\put(67.00,4.00){\makebox(0,0)[cc]{${\cdots}$}}
\put(39.00,4.00){\circle{2.00}}
\put(58.00,4.00){\line(1,0){4.00}}
\put(76.00,4.00){\line(-1,0){4.00}}
\put(77.00,4.00){\circle*{2.00}}
\put(87.00,4.00){\makebox(0,0)[cc]{${\cdots}$}}
\put(78.00,4.00){\line(1,0){4.00}}
\put(22.00,4.00){\circle{2.00}}
\put(57.00,4.00){\circle{2.00}}
\put(2.00,6.50){\makebox(0,0)[cb]{$a$}}
\put(57.00,6.50){\makebox(0,0)[cb]{$a$}}
\put(77.00,6.50){\makebox(0,0)[cb]{$a$}}
\put(31.00,2.00){\makebox(0,0)[ct]{$d$}}
\multiput(38.00,4.00)(-0.12,0.12){17}{\line(0,1){0.12}}
\multiput(38.00,4.00)(-0.12,-0.12){17}{\line(0,-1){0.12}}
\put(39.50,5.00){\makebox(0,0)[cb]{$(p+1)a/d$}}
\put(22.00,5.50){\makebox(0,0)[cb]{$pa$}}
\put(21.00,4.00){\line(-1,0){4.00}}
\put(56.00,4.00){\line(-1,0){16.00}}
\linethickness{2pt}
\put(38.00,4.00){\line(-1,0){15.00}}
\end{picture}
\end{center}

\item If the Dynkin diagram has the branching, then the nodes
at the branches not containing $\alpha_i$ are labeled by
$a,2a,\dots,pa$ and $b,2b,\dots,qb$ as above, and the nodes
at the third branch from the singularity up to $\alpha_i$ are
labeled by a decreasing arithmetic progression
$pa=qb,a+b,\dots$.
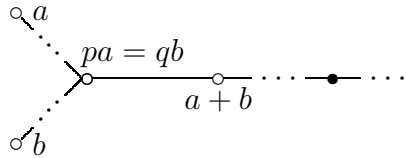
\begin{center}
\unitlength 0.40ex
\linethickness{0.4pt}
\begin{picture}(72.00,28.00)
\put(15.00,14.00){\circle{2.00}}
\put(60.00,14.00){\circle*{2.00}}
\put(50.00,14.00){\makebox(0,0)[cc]{${\cdots}$}}
\put(59.00,14.00){\line(-1,0){4.00}}
\put(15.00,14.00){\circle{2.00}}
\put(40.00,14.00){\line(1,0){5.00}}
\put(9.50,18.50){\makebox(0,0)[cc]{${\cdot}$}}
\put(8.00,20.00){\makebox(0,0)[cc]{${\cdot}$}}
\put(6.50,21.50){\makebox(0,0)[cc]{${\cdot}$}}
\multiput(14.00,14.00)(-0.12,0.12){26}{\line(0,1){0.12}}
\multiput(14.00,14.00)(-0.12,-0.12){26}{\line(0,-1){0.12}}
\put(9.50,9.50){\makebox(0,0)[cc]{${\cdot}$}}
\put(8.00,8.00){\makebox(0,0)[cc]{${\cdot}$}}
\put(6.50,6.50){\makebox(0,0)[cc]{${\cdot}$}}
\put(2.00,26.00){\circle{2.00}}
\put(2.00,2.00){\circle{2.00}}
\put(70.00,14.00){\makebox(0,0)[cc]{${\cdots}$}}
\put(61.00,14.00){\line(1,0){4.00}}
\put(5.00,26.00){\makebox(0,0)[lc]{$a$}}
\put(5.00,2.00){\makebox(0,0)[lc]{$b$}}
\put(14.00,16.00){\makebox(0,0)[lb]{$pa=qb$}}
\put(39.00,14.00){\circle{2.00}}
\multiput(5.00,5.00)(-0.11,-0.11){19}{\line(0,-1){0.11}}
\multiput(5.00,23.00)(-0.11,0.11){19}{\line(0,1){0.11}}
\put(39.00,12.50){\makebox(0,0)[ct]{$a+b$}}
\put(16.00,14.00){\line(1,0){22.00}}
\end{picture}
\end{center}

\end{itemize}

\fussy

\begin{theorem}\label{tangent}
Suppose that $G$ is simple simply connected and $P\ne G$. Then
$T_0\CE(G/\Ru{P})$ is the $(G\times L)$-stable subspace of
\begin{equation*}
\bigoplus_{i=1}^l\Hom\bigl(V(\omega_i)^{\Ru{P}},V(\omega_i)\bigr)
\end{equation*}
obtained by removing certain summands via the following
procedure:

\begin{enumerate}
\renewcommand{\theenumi}{\textup{(\arabic{enumi})}}
\renewcommand{\labelenumi}{\theenumi}

\item Take any $\alpha_k\in\Pi_L$ represented by an extreme node
of the Dynkin diagram of~$G$.

\item\label{remove} Remove subsequently all the $i$-th summands
corresponding to $\alpha_i$ which follow after $\alpha_k$ at the
Dynkin diagram until you pass the 1-st instance of
$\alpha_i\notin\Pi_L$ or the singularity.
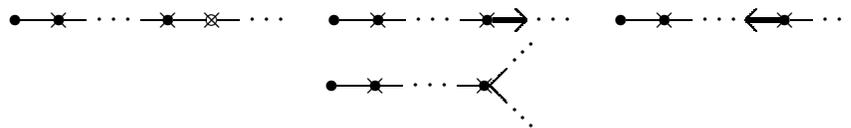
\begin{center}
\unitlength 0.40ex
\linethickness{0.4pt}
\begin{picture}(51.00,4.00)
\put(2.00,2.00){\circle*{2.00}}
\put(10.00,2.00){\circle*{2.00}}
\put(20.00,2.00){\makebox(0,0)[cc]{${\cdots}$}}
\put(38.00,2.00){\circle{2.00}}
\put(38.00,2.00){\makebox(0,0)[cc]{${\times}$}}
\put(10.00,2.00){\makebox(0,0)[cc]{${\times}$}}
\put(3.00,2.00){\line(1,0){6.00}}
\put(11.00,2.00){\line(1,0){4.00}}
\put(30.00,2.00){\circle*{2.00}}
\put(30.00,2.00){\makebox(0,0)[cc]{${\times}$}}
\put(29.00,2.00){\line(-1,0){4.00}}
\put(31.00,2.00){\line(1,0){6.00}}
\put(48.00,2.00){\makebox(0,0)[cc]{${\cdots}$}}
\put(39.00,2.00){\line(1,0){4.00}}
\end{picture}
\hfill
\unitlength 0.40ex
\linethickness{0.4pt}
\begin{picture}(45.00,4.00)
\put(2.00,2.00){\circle*{2.00}}
\put(10.00,2.00){\circle*{2.00}}
\put(20.00,2.00){\makebox(0,0)[cc]{${\cdots}$}}
\put(10.00,2.00){\makebox(0,0)[cc]{${\times}$}}
\put(3.00,2.00){\line(1,0){6.00}}
\put(11.00,2.00){\line(1,0){4.00}}
\put(30.00,2.00){\circle*{2.00}}
\put(30.00,2.00){\makebox(0,0)[cc]{${\times}$}}
\put(29.00,2.00){\line(-1,0){4.00}}
\put(42.00,2.00){\makebox(0,0)[cc]{${\cdots}$}}
\multiput(37.17,2.00)(-0.12,0.12){17}{\line(0,1){0.12}}
\multiput(37.17,2.00)(-0.12,-0.12){17}{\line(0,-1){0.12}}
\linethickness{2pt}
\put(31.00,2.00){\line(1,0){6.00}}
\end{picture}
\hfill
\unitlength 0.40ex
\linethickness{0.4pt}
\begin{picture}(45.00,4.00)
\put(2.00,2.00){\circle*{2.00}}
\put(10.00,2.00){\circle*{2.00}}
\put(20.00,2.00){\makebox(0,0)[cc]{${\cdots}$}}
\put(10.00,2.00){\makebox(0,0)[cc]{${\times}$}}
\put(3.00,2.00){\line(1,0){6.00}}
\put(11.00,2.00){\line(1,0){4.00}}
\put(42.00,2.00){\makebox(0,0)[cc]{${\cdots}$}}
\put(32.00,2.00){\circle*{2.00}}
\put(32.00,2.00){\makebox(0,0)[cc]{${\times}$}}
\put(33.00,2.00){\line(1,0){4.00}}
\multiput(24.83,2.00)(0.12,0.12){17}{\line(0,1){0.12}}
\multiput(24.83,2.00)(0.12,-0.12){17}{\line(0,-1){0.12}}
\linethickness{2pt}
\put(31.00,2.00){\line(-1,0){6.00}}
\end{picture}
\\
\unitlength 0.40ex
\linethickness{0.4pt}
\begin{picture}(40.00,18.50)
\put(2.00,9.00){\circle*{2.00}}
\put(10.00,9.00){\circle*{2.00}}
\put(20.00,9.00){\makebox(0,0)[cc]{${\cdots}$}}
\put(10.00,9.00){\makebox(0,0)[cc]{${\times}$}}
\put(3.00,9.00){\line(1,0){6.00}}
\put(11.00,9.00){\line(1,0){4.00}}
\put(30.00,9.00){\circle*{2.00}}
\put(30.00,9.00){\makebox(0,0)[cc]{${\times}$}}
\put(29.00,9.00){\line(-1,0){4.00}}
\put(35.50,13.50){\makebox(0,0)[cc]{${\cdot}$}}
\put(37.00,15.00){\makebox(0,0)[cc]{${\cdot}$}}
\put(38.50,16.50){\makebox(0,0)[cc]{${\cdot}$}}
\multiput(31.00,9.00)(0.12,0.12){26}{\line(0,1){0.12}}
\multiput(31.00,9.00)(0.12,-0.12){26}{\line(0,-1){0.12}}
\put(35.50,4.50){\makebox(0,0)[cc]{${\cdot}$}}
\put(37.00,3.00){\makebox(0,0)[cc]{${\cdot}$}}
\put(38.50,1.50){\makebox(0,0)[cc]{${\cdot}$}}
\end{picture}
\end{center}

\item If $G$ is simply laced, and at least two branches of the
Dynkin diagram are contained in~$\Pi_L$, then continue removing
the summands along the 3-rd branch after the singularity as
in~\ref{remove} until, in the case $G=\Ee_l$, the removed
segment becomes longer than both other branches.
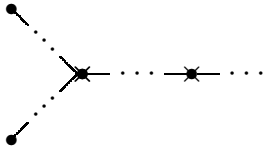
\begin{center}
\unitlength 0.40ex
\linethickness{0.4pt}
\begin{picture}(47.00,29.00)
\put(15.00,14.00){\circle*{2.00}}
\put(15.00,14.00){\makebox(0,0)[cc]{${\times}$}}
\put(35.00,14.00){\circle*{2.00}}
\put(25.00,14.00){\makebox(0,0)[cc]{${\cdots}$}}
\put(35.00,14.00){\makebox(0,0)[cc]{${\times}$}}
\put(34.00,14.00){\line(-1,0){4.00}}
\put(15.00,14.00){\circle*{2.00}}
\put(15.00,14.00){\makebox(0,0)[cc]{${\times}$}}
\put(16.00,14.00){\line(1,0){4.00}}
\put(9.50,18.50){\makebox(0,0)[cc]{${\cdot}$}}
\put(8.00,20.00){\makebox(0,0)[cc]{${\cdot}$}}
\put(6.50,21.50){\makebox(0,0)[cc]{${\cdot}$}}
\multiput(14.00,14.00)(-0.12,0.12){26}{\line(0,1){0.12}}
\multiput(14.00,14.00)(-0.12,-0.12){26}{\line(0,-1){0.12}}
\put(9.50,9.50){\makebox(0,0)[cc]{${\cdot}$}}
\put(8.00,8.00){\makebox(0,0)[cc]{${\cdot}$}}
\put(6.50,6.50){\makebox(0,0)[cc]{${\cdot}$}}
\put(2.00,26.00){\circle*{2.00}}
\multiput(2.00,26.00)(0.12,-0.12){26}{\line(0,-1){0.12}}
\put(2.00,2.00){\circle*{2.00}}
\multiput(2.00,2.00)(0.12,0.12){26}{\line(0,1){0.12}}
\put(45.00,14.00){\makebox(0,0)[cc]{${\cdots}$}}
\put(36.00,14.00){\line(1,0){4.00}}
\end{picture}
\end{center}

\item If $G$ is not simply laced, and you have passed the
singularity along the direction to long roots, then continue
removing summands as in~\ref{remove}.
\begin{center}
\unitlength 0.40ex
\linethickness{0.4pt}
\begin{picture}(64.00,6.00)
\put(2.00,2.00){\circle*{2.00}}
\put(12.00,2.00){\makebox(0,0)[cc]{${\cdots}$}}
\put(3.00,2.00){\line(1,0){4.00}}
\put(34.00,2.00){\makebox(0,0)[cc]{${\cdots}$}}
\put(24.00,2.00){\circle*{2.00}}
\put(24.00,2.00){\makebox(0,0)[cc]{${\times}$}}
\put(25.00,2.00){\line(1,0){4.00}}
\multiput(16.83,2.00)(0.12,0.12){17}{\line(0,1){0.12}}
\multiput(16.83,2.00)(0.12,-0.12){17}{\line(0,-1){0.12}}
\put(44.00,2.00){\circle*{2.00}}
\put(44.00,2.00){\makebox(0,0)[cc]{${\times}$}}
\put(43.00,2.00){\line(-1,0){4.00}}
\put(52.00,2.00){\circle{2.00}}
\put(52.00,2.00){\makebox(0,0)[cc]{${\times}$}}
\put(62.00,2.00){\makebox(0,0)[cc]{${\cdots}$}}
\put(45.00,2.00){\line(1,0){6.00}}
\put(53.00,2.00){\line(1,0){4.00}}
\linethickness{2pt}
\put(23.00,2.00){\line(-1,0){6.00}}
\end{picture}
\end{center}

\end{enumerate}

\end{theorem}

\begin{Examples}
Let $G=\Ee_8$ and $P$ be the projective stabilizer of a
highest weight vector in~$V(\omega_1)$, in the enumeration
of~\cite[Table~1]{sem}. Then $L'=\Ee_7$, with the simple roots
corresponding to the black nodes of the diagram:
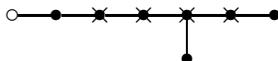
\begin{center}
\unitlength 0.40ex
\linethickness{0.4pt}
\begin{picture}(51.00,11.00)
\put(2.00,10.00){\circle{2.00}}
\put(10.00,10.00){\circle*{2.00}}
\put(3.00,10.00){\line(1,0){6.00}}
\put(18.00,10.00){\circle*{2.00}}
\put(26.00,10.00){\circle*{2.00}}
\put(34.00,10.00){\circle*{2.00}}
\put(42.00,10.00){\circle*{2.00}}
\put(50.00,10.00){\circle*{2.00}}
\put(18.00,10.00){\makebox(0,0)[cc]{${\times}$}}
\put(26.00,10.00){\makebox(0,0)[cc]{${\times}$}}
\put(34.00,10.00){\makebox(0,0)[cc]{${\times}$}}
\put(42.00,10.00){\makebox(0,0)[cc]{${\times}$}}
\put(11.00,10.00){\line(1,0){6.00}}
\put(19.00,10.00){\line(1,0){6.00}}
\put(27.00,10.00){\line(1,0){6.00}}
\put(35.00,10.00){\line(1,0){6.00}}
\put(43.00,10.00){\line(1,0){6.00}}
\put(34.00,2.00){\circle*{2.00}}
\put(34.00,3.00){\line(0,1){6.00}}
\end{picture}
\end{center}
We have $\dim\CE(G/\Ru{P})=191$, but the minimal ambient
$G$-module is
\begin{multline*}
\bigl(V(\omega_1)\otimes V_L(\omega_1)^{*}\bigr)
\oplus\bigl(V(\omega_2)\otimes V_L(\omega_2)^{*}\bigr)\\
\oplus\bigl(V(\omega_7)\otimes V_L(\omega_7)^{*}\bigr)
\oplus\bigl(V(\omega_8)\otimes V_L(\omega_8)^{*}\bigr)
\end{multline*}
of dimension
$248\cdot1+30380\cdot56+3875\cdot133+147250\cdot912=136508903$.

Now take $G=\Ff_4$ and $P$ the projective stabilizer of a
highest weight vector in~$V(\omega_1)$ again. Then
$L'=\Spin(7)$, with the simple roots corresponding to the black
nodes of the diagram:
\begin{center}
\unitlength 0.40ex
\linethickness{0.4pt}
\begin{picture}(27.00,3.67)
\put(2.00,2.00){\circle{2.00}}
\put(10.00,2.00){\circle*{2.00}}
\put(3.00,2.00){\line(1,0){6.00}}
\put(18.00,2.00){\circle*{2.00}}
\put(26.00,2.00){\circle*{2.00}}
\put(18.00,2.00){\makebox(0,0)[cc]{${\times}$}}
\put(19.00,2.00){\line(1,0){6.00}}
\put(10.50,2.50){\line(1,0){7.00}}
\put(10.50,1.50){\line(1,0){7.00}}
\multiput(11.00,2.00)(0.12,0.12){14}{\line(1,0){0.12}}
\multiput(11.00,2.00)(0.12,-0.12){14}{\line(1,0){0.12}}
\end{picture}
\end{center}
We have $\dim\CE(G/\Ru{P})=37$, and the minimal ambient
$G$-module
\begin{equation*}
\bigl(V(\omega_1)\otimes V_L(\omega_1)^{*}\bigr)\oplus
\bigl(V(\omega_2)\otimes V_L(\omega_2)^{*}\bigr)\oplus
\bigl(V(\omega_4)\otimes V_L(\omega_4)^{*}\bigr)
\end{equation*}
has dimension
$26\cdot1+273\cdot8+52\cdot7=2574$.
\end{Examples}

\begin{proof}
The space
$\bigoplus_i\Hom\bigl(V(\omega_i)^{\Ru{P}},V(\omega_i)\bigr)$ is
dual to $\bigoplus_i\kk[G/\Ru{P}]_{(\omega_i)}$, a generating
subspace of $\kk[G/\Ru{P}]$. To obtain the tangent space, it
suffices to remove summands corresponding to $\omega_i$ which
are $L$-generated by the others.

First observe that if $\alpha_i\notin\Pi_L\sqcup\partial\Pi_L$,
then $\omega_i$ is not $L$-generated by the other fundamental
weights. Indeed, specialize the multi-grading of $\kk[G/\Ru{P}]$
using $\dual\alpha_i$. Then $\deg\omega_i=1$, but
$\deg\omega_j=0$, $\forall j\ne i$.

Secondly, $\omega_i$~is $L$-generated by the other $\omega_j$'s
iff it is $L_k$-generated by the other $\omega_j$'s such that
$\alpha_j\in\Pi_{L_k}$, where $L_k$ is one of the simple
factors of~$L$. Indeed, each dominant weight $L$-generated by
$\omega_j$'s is the sum of dominant weights $L_k$-generated by
$\omega_j$ such that $\alpha_j\in\Pi_{L_k}$, over all simple
factors $L_k\subseteq L$, and of a dominant weight generated by
$\omega_j$ such that $\alpha_j\notin\setminus\Pi_L$. However,
specializing the multi-grading of $\kk[G/\Ru{P}]$ to a
non-negative grading such that $\kk[G/\Ru{P}]_0=\kk$
(Remark~\ref{grading}) shows that $\omega_j$'s do not
$L_k$-generate~$0$. The assertion follows, because $\omega_i$
cannot be decomposed as a non-trivial sum of dominant weights.
Thus we may assume that $\Pi_L$ is indecomposable.

In order to verify that certain $\omega_j$ are $L$-generated by
the others (as asserted in Theorem~\ref{tangent}), we use the
following formul{\ae} \cite[Table~5]{sem}:
\begin{enumerate}
\renewcommand{\theenumi}{\textup{(\arabic{enumi})}}
\renewcommand{\labelenumi}{\theenumi}

\item\label{1^i} $V(\resomega_1)^{\otimes i}\hookleftarrow
V(\resomega_i)$ for:\\ $L'=\SL(m),\Sp(2m)$, $1\le i\le m$;
$L'=\Spin(2m+1)$, $1\le i<m$; $L'=\Spin(2m)$, $1\le i\le m-2$.

\item\label{1*(1^*)} $V(\resomega_1)\otimes V(\resomega_{m-1})
\hookleftarrow V(0)$ for $L'=\SL(m)$.

\item\label{1^2} $V(\resomega_1)^{\otimes 2}\hookleftarrow V(0)$
for $L'=\Sp(2m),\Spin(2m)$.

\item\label{1*(m-1)} $V(\resomega_1)\otimes V(\resomega_{m-1})
\hookleftarrow V(\resomega_m)$ for $L'=\Spin(2m)$.

\item\label{spin^2} $V(\resomega_m)^{\otimes 2}=
V(2\resomega_m)\oplus V(\resomega_{m-1})\oplus\dots\oplus
V(\resomega_1)\oplus V(0)$\\ for $L'=\Spin(2m+1)$.

\item\label{(spin+)^2} $V(\resomega_m)^{\otimes 2}=
V(2\resomega_m)\oplus V(\resomega_{m-2})\oplus
V(\resomega_{m-4})\oplus\cdots$\\ for $L'=\Spin(2m)$.

\item\label{(spin-)*(spin+)} $V(\resomega_{m-1})\otimes
V(\resomega_m)=V(\resomega_{m-1}+\resomega_m)\oplus
V(\resomega_{m-3})\oplus V(\resomega_{m-5})\oplus\cdots$\\ for
$L'=\Spin(2m)$.

\item\label{(m-1)^2} $V(\resomega_{m-1})^{\otimes 2}
\hookleftarrow V(\resomega_{m-2})\oplus V(\resomega_{m-5})$ for
$L'=\Ee_m$, $m=6,7$.

\item\label{m^2} $V(\resomega_m)^{\otimes 2}
\hookleftarrow V(\resomega_{m-3})$ for $L'=\Ee_m$,
$m=6,7$.

\item\label{(m-1)*m} $V(\resomega_{m-1})\otimes V(\resomega_m)
\hookleftarrow V(\resomega_{m-4})$ for $L'=\Ee_m$,
$m=6,7$.

\end{enumerate}
Here the fundamental weights of $L$ are numbered according to
\cite[Table~1]{sem}. The respective relations between degrees
are easily verified using the above description of degrees
w.r.t.\ fundamental coweights. Note that it suffices to
consider degrees w.r.t.\ $\dual\omega_i$ such that
$\alpha_i\in\partial\Pi_L$, because
$\X_{*}(Z)\otimes\QQ=\langle\dual\omega_i,\dual\alpha_j\mid
\alpha_i\in\partial\Pi_L,\
\alpha_j\notin\Pi_L\sqcup\partial\Pi_L\rangle$ and the degrees
of fundamental weights corresponding to roots in
$\Pi_L\sqcup\partial\Pi_L$ w.r.t.\
$\dual\alpha_j\notin\dual\Pi_L\sqcup\partial\dual\Pi_L$ are
zero.

For instance, suppose that $\Pi=\Ee_l$,
$\Pi_L=\Dd_{l-1}$. Let us enumerate the simple roots of
$G$ as at the picture:
\begin{center}
\unitlength 0.40ex
\linethickness{0.4pt}
\begin{picture}(80.00,17.00)
\put(2.00,13.00){\circle*{2.00}}
\put(10.00,13.00){\circle*{2.00}}
\put(20.00,13.00){\makebox(0,0)[cc]{${\cdots}$}}
\put(78.00,13.00){\circle{2.00}}
\put(78.00,13.00){\makebox(0,0)[cc]{${\times}$}}
\put(10.00,13.00){\makebox(0,0)[cc]{${\times}$}}
\put(3.00,13.00){\line(1,0){6.00}}
\put(11.00,13.00){\line(1,0){4.00}}
\put(43.00,13.00){\circle*{2.00}}
\put(43.00,13.00){\makebox(0,0)[cc]{${\times}$}}
\put(43.00,2.00){\circle*{2.00}}
\put(66.00,13.00){\circle*{2.00}}
\put(66.00,13.00){\makebox(0,0)[cc]{${\times}$}}
\put(2.00,11.00){\makebox(0,0)[ct]{$\alpha_1$}}
\put(2.00,15.50){\makebox(0,0)[cb]{$a$}}
\put(10.00,11.00){\makebox(0,0)[ct]{$\alpha_2$}}
\put(10.00,15.50){\makebox(0,0)[cb]{$2a$}}
\put(42.00,11.00){\makebox(0,0)[rt]{$\alpha_{l-3}$}}
\put(41.00,14.00){\makebox(0,0)[cb]{$(l-3)a=2b$}}
\put(66.00,11.00){\makebox(0,0)[ct]{$\alpha_{l-1}$}}
\put(66.00,15.00){\makebox(0,0)[cb]{$a+b$}}
\put(78.00,11.00){\makebox(0,0)[ct]{$\alpha_l$}}
\put(78.00,15.50){\makebox(0,0)[cb]{$2a$}}
\put(41.00,2.00){\makebox(0,0)[rc]{$b$}}
\put(45.00,2.00){\makebox(0,0)[lc]{$\alpha_{l-2}$}}
\put(43.00,3.00){\line(0,1){9.00}}
\put(77.00,13.00){\line(-1,0){10.00}}
\put(65.00,13.00){\line(-1,0){21.00}}
\put(42.00,13.00){\line(-1,0){17.00}}
\end{picture}
\end{center}
We consider the degrees w.r.t.~$\dual\omega_l$. Using \ref{1^i}
and $d_i=id_1$, we verify that $\omega_i$ are $L$-generated
by~$\omega_1$, $1\le i\le l-3$. By \ref{1*(m-1)} and
$d_1+d_{l-2}=d_{l-1}$, we see that $\omega_{l-1}$ is
$L$-generated by $\omega_1,\omega_{l-2}$. Finally, \ref{1^2}~and
$d_l=2d_1$ implies that $\omega_l$ is $L$-generated
by~$\omega_1$.

It remains to prove that the remaining fundamental weights are
not $L$-generated by the others. We shall use
the following observation from the
representation theory of $\SL(m)$:
\begin{gather*}
\tag{\degs}
V(\resomega_{j_1})\otimes\dots\otimes V(\resomega_{j_n})
\hookleftarrow V(\resomega_i)\implies j_1+\dots+j_n\ge i\\
1\le i,j_1,\dots,j_n\le m
\end{gather*}
(Here $\resomega_m=0$ and the other $\resomega_j$ are the
fundamental weights of $\SL(m)$ in the standard order.) In the
sequel, we shall frequently apply (\degs) to $L'=\SL(m)$ in the
following way: it often happens that the conclusion of (\degs)
implies $d_{j_1}+\dots+d_{j_n}>d_i$, whence $\omega_i$ is not
$L$-generated by the other $\omega_j$'s.

First suppose that the Dynkin diagram of $G$ has no branching.

Fix any $\alpha_m\notin\Pi_L$ and consider the degrees of
fundamental weights w.r.t.\ $\dual\omega_m$ on one of the
segments of $\Pi\setminus\{\alpha_m\}$. From the above
description of degrees, we easily see that $d_i<d_j+d_k$
unless $\alpha_j,\alpha_k$ are further from $\alpha_m$
than~$\alpha_i$. Hence on this segment each $\omega_i$ could be
$L$-generated only by fundamental weights corresponding to roots
on the other side from $\alpha_i$ than~$\alpha_m$. We
immediately deduce that if $\Pi_L$ does not contain an extreme
node of (the Dynkin diagram of) $\Pi$, then no fundamental
weights are $L$-generated by the others.

Now assume that $\alpha_m$ is a short root and look at the
degrees on the segment from an extreme node to $\alpha_m$
containing the singularity. For
$\Pi=\Cc_l,\Ff_4,\Gg_2$ we have $d_i<d_j+d_k$ whenever
$\alpha_i$ is short, hence fundamental weights corresponding to
short roots after the singularity up to $\partial\Pi_L$ are not
$L$-generated by the others. The same assertion for the unique
short root $\alpha_l$ of $\Pi=\Bb_l$ stems
from~(\degs).

Next, suppose that the Dynkin diagram of $G$ has the
branching. We consider the degrees w.r.t.\ $\dual\omega_m$
such that $\alpha_m$ corresponds to the extreme node of a ray of
the Dynkin diagram. For convenience of the reader, let us
indicate these degrees at the diagrams, where the black node
corresponds to~$\alpha_m$ (the picture for $\alpha_m$ on the
long ray of $\Ee_l$ is obtained from that for
$\Ee_8$ by cutting off $8-l$ subsequent nodes on the long
ray, starting with the extreme node):
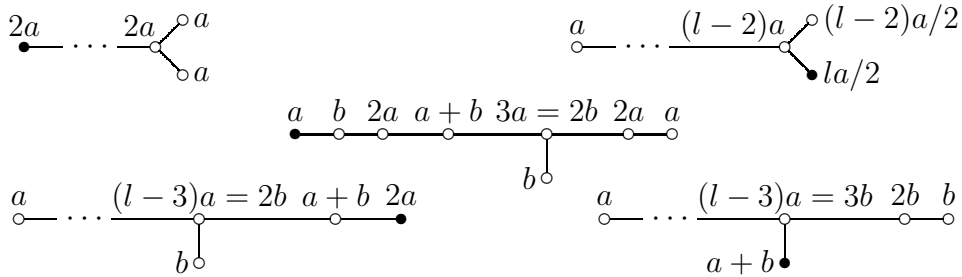
\begin{center}
\unitlength 0.40ex
\linethickness{0.4pt}
\begin{picture}(35.00,14.00)
\put(3.00,7.00){\circle*{2.00}}
\put(15.00,7.00){\makebox(0,0)[cc]{${\cdots}$}}
\put(27.00,7.00){\circle{2.00}}
\put(27.00,9.00){\makebox(0,0)[rb]{$2a$}}
\put(3.00,9.00){\makebox(0,0)[cb]{$2a$}}
\put(32.00,12.00){\circle{2.00}}
\put(34.00,12.00){\makebox(0,0)[lc]{$a$}}
\put(32.00,2.00){\circle{2.00}}
\put(34.00,2.00){\makebox(0,0)[lc]{$a$}}
\multiput(31.33,11.33)(-0.12,-0.12){31}{\line(0,-1){0.12}}
\multiput(31.33,2.67)(-0.12,0.12){31}{\line(0,1){0.12}}
\put(4.00,7.00){\line(1,0){6.00}}
\put(26.00,7.00){\line(-1,0){6.00}}
\end{picture}
\hfill
\unitlength 0.40ex
\linethickness{0.4pt}
\begin{picture}(72.00,14.00)
\put(2.00,7.00){\circle{2.00}}
\put(14.00,7.00){\makebox(0,0)[cc]{${\cdots}$}}
\put(40.00,7.00){\circle{2.00}}
\put(40.00,8.00){\makebox(0,0)[rb]{$(l-2)a$}}
\put(2.00,9.50){\makebox(0,0)[cb]{$a$}}
\put(45.00,12.00){\circle{2.00}}
\put(47.00,12.00){\makebox(0,0)[lc]{$(l-2)a/2$}}
\put(45.00,2.00){\circle*{2.00}}
\put(47.00,2.00){\makebox(0,0)[lc]{$la/2$}}
\multiput(44.33,11.33)(-0.12,-0.12){31}{\line(-1,0){0.12}}
\multiput(44.33,2.67)(-0.12,0.12){31}{\line(-1,0){0.12}}
\put(3.00,7.00){\line(1,0){6.00}}
\put(39.00,7.00){\line(-1,0){20.00}}
\end{picture}
\unitlength 0.40ex
\linethickness{0.4pt}
\begin{picture}(73.00,18.50)
\put(2.00,10.00){\circle*{2.00}}
\put(3.00,10.00){\line(1,0){6.00}}
\put(10.00,10.00){\circle{2.00}}
\put(18.00,10.00){\circle{2.00}}
\put(30.00,10.00){\circle{2.00}}
\put(48.00,10.00){\circle{2.00}}
\put(63.00,10.00){\circle{2.00}}
\put(71.00,10.00){\circle{2.00}}
\put(11.00,10.00){\line(1,0){6.00}}
\put(64.00,10.00){\line(1,0){6.00}}
\put(48.00,2.00){\circle{2.00}}
\put(48.00,3.00){\line(0,1){6.00}}
\put(2.00,12.50){\makebox(0,0)[cb]{$a$}}
\put(10.00,12.50){\makebox(0,0)[cb]{$b$}}
\put(18.00,12.50){\makebox(0,0)[cb]{$2a$}}
\put(30.00,12.00){\makebox(0,0)[cb]{$a+b$}}
\put(48.00,12.50){\makebox(0,0)[cb]{$3a=2b$}}
\put(63.00,12.50){\makebox(0,0)[cb]{$2a$}}
\put(71.00,12.50){\makebox(0,0)[cb]{$a$}}
\put(46.00,2.00){\makebox(0,0)[rc]{$b$}}
\put(19.00,10.00){\line(1,0){10.00}}
\put(31.00,10.00){\line(1,0){16.00}}
\put(49.00,10.00){\line(1,0){13.00}}
\end{picture}
\\
\unitlength 0.40ex
\linethickness{0.4pt}
\begin{picture}(74.00,15.00)
\put(72.00,10.00){\circle*{2.00}}
\put(60.00,10.00){\circle{2.00}}
\put(35.00,10.00){\circle{2.00}}
\put(2.00,10.00){\circle{2.00}}
\put(35.00,2.00){\circle{2.00}}
\put(35.00,3.00){\line(0,1){6.00}}
\put(72.00,12.50){\makebox(0,0)[cb]{$2a$}}
\put(60.00,12.00){\makebox(0,0)[cb]{$a+b$}}
\put(35.00,11.00){\makebox(0,0)[cb]{$(l-3)a=2b$}}
\put(2.00,12.50){\makebox(0,0)[cb]{$a$}}
\put(71.00,10.00){\line(-1,0){10.00}}
\put(33.00,2.00){\makebox(0,0)[rc]{$b$}}
\put(36.00,10.00){\line(1,0){23.00}}
\put(34.00,10.00){\line(-1,0){15.00}}
\put(14.00,10.00){\makebox(0,0)[cc]{${\cdots}$}}
\put(3.00,10.00){\line(1,0){6.00}}
\end{picture}
\hfill
\unitlength 0.40ex
\linethickness{0.4pt}
\begin{picture}(67.00,15.00)
\put(65.00,10.00){\circle{2.00}}
\put(57.00,10.00){\circle{2.00}}
\put(35.00,10.00){\circle{2.00}}
\put(2.00,10.00){\circle{2.00}}
\put(35.00,2.00){\circle*{2.00}}
\put(35.00,3.00){\line(0,1){6.00}}
\put(65.00,12.50){\makebox(0,0)[cb]{$b$}}
\put(57.00,12.50){\makebox(0,0)[cb]{$2b$}}
\put(35.00,11.00){\makebox(0,0)[cb]{$(l-3)a=3b$}}
\put(2.00,12.50){\makebox(0,0)[cb]{$a$}}
\put(33.00,2.00){\makebox(0,0)[rc]{$a+b$}}
\put(34.00,10.00){\line(-1,0){15.00}}
\put(14.00,10.00){\makebox(0,0)[cc]{${\cdots}$}}
\put(3.00,10.00){\line(1,0){6.00}}
\put(64.00,10.00){\line(-1,0){6.00}}
\put(56.00,10.00){\line(-1,0){20.00}}
\end{picture}
\end{center}
We enumerate the simple roots and fundamental weights of $G$
according to \cite[Table~1]{sem}.

One sees from this picture that the fundamental weights
$\omega_i$ corresponding to extreme nodes of the two other rays,
as well as $\omega_i$ for $\alpha_m$ on the long ray of
$\Ee_l$ and $i<l-5$, are not $L$-generated by the
others. Indeed, for $\alpha_m$ on a short ray of $\Dd_l$
or on the middle ray of $\Ee_l$ we can use (\degs) or the
fact that a semispinor weight of $\Dd_{l-1}$ is not
$\Dd_{l-1}$-generated by~$\resomega_1$. In all other
cases we have $d_i<d_j+d_k$, $1\le j,k\le l$.

Now assume that the extreme nodes of at least two rays are not
in~$\Pi_L$.

If one of these rays is long, then looking at the degrees
w.r.t.\ the respective $\dual\omega_m$ shows that, except for
the weights on the 3-rd ray, which could be $L$-generated by the
one at the extreme node, the only possibilities for
$L$-generation are: $\Pi=\Dd_l$, $d_i=2d_l$,
$V_L(\omega_i)\hookrightarrow V_L(\omega_{l-1})^{\otimes2}$, $1\le
i<l-2$; $\Pi=\Ee_8$, $d_5=2d_2$,
$V_L(\omega_5)\hookrightarrow V_L(\omega_2)^{\otimes2}$;
$\Pi=\Ee_8$, $d_4=d_2+d_7$, $V_L(\omega_4)\hookrightarrow
V_L(\omega_2)\otimes V_L(\omega_7)$; $\Pi=\Ee_l$,
$d_{l-5}=2d_{l-1}$, $V_L(\omega_{l-5})\hookrightarrow
V_L(\omega_{l-1})^{\otimes2}$. However these possibilities are
excluded by considering the degrees w.r.t.\ the extreme node of
the 2-nd ray.

Otherwise, consider the degrees w.r.t.\ $\dual\omega_m$ such
that $\alpha_m$ is at the extreme node of the short ray.

For $\Pi=\Dd_l$, $V_L(\omega_i)\hookrightarrow
V_L(\omega_{j_1})\otimes\dots\otimes V_L(\omega_{j_n})$,
$j_1,\dots,j_n\le l-2<i$, implies $d_i=d_{j_1}+\dots+d_{j_n}$.
However, considering the degrees w.r.t.\ the extreme node of the
2-nd ray violates this equality, a contradiction.

For $\Pi=\Ee_l$, the only possibilities for
$L$-generation are: $l=6$, $d_4=d_6=2d_1$, and $V_L(\omega_4)$
or $V_L(\omega_6)$ is contained in $V_L(\omega_1)^{\otimes2}$.
However these possibilities are excluded by considering the
degrees w.r.t.\ the extreme node of the 2-nd ray.

We conclude that no fundamental weights on a segment between two
nodes of the Dynkin diagram not contained in $\Pi_L$ are
$L$-generated by the others, except possibly the one at the
singularity provided that one of the rays is contained
in~$\Pi_L$. This completes the proof.
\end{proof}

\begin{remark}
Our results immediately extend to the case, where $G$ is
semisimple simply connected, see Remark~\ref{G->simple}. The
general case looks more complicated, because the structure of
$\X^{+}$ is more involved. It would be interesting to solve the
problem in full generality.
\end{remark}



\begin{thebibliography}{PRV67}



\bibitem[Ar03]{ar2}
I. V. Arzhantsev, \emph{Algebras with finitely generated invariant
subalgebras}, Ann. Inst. Fourier {\bf53:2} (2003), 379--398.

\bibitem[AT01]{at}
I. V. Arzhantsev and D. A. Timashev, \emph{Affine embeddings with a finite
number of orbits}, Transformation Groups {\bf6:2} (2001),
101--110.


\bibitem[Br97]{var.sph}
M.~Brion, \emph{Vari{\'e}t{\'e}s sph{\'e}riques}, Notes de la
session de la S. M. F. ``Op{\'e}rations hamiltoniennes et
op{\'e}rations de groupes alg{\'e}briques'', Grenoble, 1997;\\
\texttt{http://www-fourier.ujf-grenoble.fr/\~{}mbrion/spheriques.ps}.

\bibitem[Gr73]{gr1}
F. D. Grosshans, \emph{Observable groups and Hilbert's fourteenth problem},
Amer. J. Math. {\bf95:1} (1973), 229--253.

\bibitem[Gr83]{gr2}
F. D. Grosshans, \emph{The Invariants of Unipotent Radicals of Parabolic
Subgroups}, Invent. Math. {\bf73} (1983), 1--9.

\bibitem[Gr97]{gr}
F. D. Grosshans, \emph{Algebraic Homogeneous Spaces and Invariant Theory},
Lect. Notes Math. {\bf1673}, Springer-Verlag, Berlin, 1997.

\bibitem[Ke78]{ke}
G. Kempf, \emph{Instability in invariant theory}, Ann. of Math.
{\bf108:2} (1978), 299--316.

\bibitem[Kn91]{sph}
F.~Knop, \emph{The Luna--Vust theory of spherical embeddings}, Proc. Hyderabad
Conf. on Algebraic Groups (S.~Ramanan, ed.), pp. 225--249, Manoj Prakashan,
Madras, 1991.

\bibitem[Kr85]{inv}
H.~Kraft, \emph{Geometrische Methoden in der Invariantentheorie},
Vieweg, Braunschweig--Wiesbaden, 1985.

\bibitem[Lu73]{lu2}
D. Luna, \emph{Slices {\'e}tales}, M{\'e}m. Soc. Math. France (N.S.)
{\bf33} (1973), 81--105.

\bibitem[Lu75]{lu1}
D. Luna, \emph{Adh\'erences d'orbite et invariants}, Invent. Math.
{\bf29} (1975), 231--238.

\bibitem[Ma60]{mat}
Y. Matsushima, \emph{Espaces homog\`enes de Stein des groupes de
Lie complexes}, Nagoya Math. J. {\bf16} (1960), 205--216.

\bibitem[On60]{on}
A. L. Onishchik, \emph{Complex envelopes of compact homogeneous spaces},
Dokl. Akad. Nauk SSSR {\bf130:4} (1960), 726--729 (in Russian).

\bibitem[On62]{Aut(C-flag)}
A. L. Onishchik, \emph{Inclusion relations between transitive
compact transformation groups}, Trudy Moskov. Mat. Obshch.
{\bf11} (1962), 199--242 (in Russian).

\bibitem[OV88]{sem}
A.~L. Onishchik, E.~B. Vinberg, \emph{Seminar on Lie groups and
algebraic groups}, Nauka, Moscow, 1988 (in Russian);
English transl.: \emph{Lie groups and algebraic groups},
Springer-Verlag, Berlin--Heidelberg--New York, 1990.

\bibitem[PRV67]{PRV}
K.~R.~Parthasarathy, R.~Ranga Rao, V.~S.~Varadarajan,
\emph{Representations of complex semi-simple Lie groups and Lie
algebras}, Ann. of Math. {\bf85} (1967), 383--429.

\bibitem[Po73]{pop}
V.~L.~Popov, \emph{Quasihomogeneous affine algebraic varieties of the
group $\SL_2$}, Izv. Akad. Nauk SSSR, Ser. Mat.
{\bf37} (1973), 792--832 (in Russian);
English transl.: Math. USSR-Izv. {\bf7} (1973), 793--831.

\bibitem[PV89]{pv2}
V. L. Popov and E. B. Vinberg, \emph{Invariant Theory}, Itogy Nauki i
Tekhniki, Sovr. Problemy Mat., Fund. Napravlenia, vol. 55,
VINITI, Moscow, 1989, pp. 137--309 (in Russian); English
transl.: Algebraic Geometry IV, Encyclopaedia of Math.
Sciences, vol. 55, Springer-Verlag, Berlin, 1994, pp. 123--278.

\bibitem[Ri77]{rich}
R. W. Richardson, \emph{Affine coset spaces of reductive algebraic groups},
Bull. London Math. Soc. {\bf9} (1977), 38--41.

\bibitem[Rit98]{mon}
A.~Rittatore, \emph{Algebraic monoids and group embeddings},
Transformation Groups {\bf 3:4} (1998), 375--396.

\bibitem[St82]{Aut(flag)}
M. Steinsiek, \emph{Transformation groups on
homogeneous-rational manifolds}, Math. Ann. {\bf260:4} (1982),
423--435.

\bibitem[Su88]{sukh}
A.~A.~Sukhanov, \emph{Description of the observable subgroups of
linear algebraic groups},
Mat. Sbornik {\bf137:1} (1988), 90--102 (in Russian);
English transl.: Math. USSR-Sb. {\bf65:1} (1990), 97--108.

\bibitem[Ti03]{grp.comp}
D.~A. Timashev, \emph{Equivariant compactifications of reductive groups},
Mat. Sbornik {\bf194:4} (2003), 119--146 (in Russian);
English transl.: Sbornik: Mathematics {\bf194:4} (2003), 589--616;
\texttt{arXiv:math.AG/0207034}.

\bibitem[Ti03$'$]{eq.emb}
D.~A. Timashev, \emph{Equivariant embeddings of homogeneous spaces},
Pr{\'e}\-pub\-li\-ca\-tion de l'Institut Fourier
n$^{\circ}$~611, 2003.

\bibitem[Vi95]{vi}
E.~B.~Vinberg, \emph{On reductive algebraic semigroups}, Lie Groups and
Lie Algebras: E.~B.~Dynkin Seminar (S.~Gindikin, E.~Vinberg, eds.),
AMS Transl. {\bf169} (1995), 145--182.

\end{thebibliography}
\end{document}